\documentclass[12pt,a4paper]{article}

\usepackage[british]{babel}

\usepackage[T1]{fontenc}

\usepackage[utf8]{inputenc}

\usepackage[sc]{mathpazo}

\usepackage{csquotes}

\usepackage[backend=biber, maxbibnames=10, maxcitenames=4, maxalphanames=4, bibencoding=utf8, style=alphabetic, isbn=false, url=false, doi=true]{biblatex}
\renewbibmacro{in:}{}
\renewbibmacro*{doi+eprint+url}{%
	\printfield{doi}%
	\newunit\newblock%
	\iftoggle{bbx:eprint}{%
		\usebibmacro{eprint}%
	}{}%
	\newunit\newblock%
	\iffieldundef{doi}{%
		\usebibmacro{url+urldate}}%
	{}%
}

\usepackage{setspace}

\usepackage[indentafter]{titlesec}
\titleformat{name=\section}{}{\thetitle.}{0.8em}{\centering\scshape}
\titleformat{name=\subsection}[runin]{}{\thetitle.}{0.5em}{\bfseries}[.]
\titleformat{name=\subsubsection}[runin]{}{\thetitle.}{0.5em}{\itshape}[.]
\titleformat{name=\paragraph,numberless}[runin]{}{}{0em}{}[.]
\titlespacing{\paragraph}{0em}{0em}{0.5em}
\titleformat{name=\subparagraph,numberless}[runin]{}{}{0em}{}[.]
\titlespacing{\subparagraph}{0em}{0em}{0.5em}

\usepackage{url}

\usepackage{mathtools}

\usepackage{amssymb}

\usepackage{amsthm}

\usepackage{mathrsfs}

\usepackage{dsfont}

\usepackage{enumitem}
\setlist[enumerate]{noitemsep, partopsep=0pt, topsep=0pt, parsep=0pt, itemsep=0pt}
\setlist[itemize]{noitemsep, partopsep=0pt, topsep=0pt, parsep=0pt, itemsep=0pt}

\usepackage{interval}

\intervalconfig{soft open fences}

\usepackage[normalem]{ulem}

\usepackage[textsize=tiny, obeyDraft]{todonotes}

\setuptodonotes{fancyline}
\usepackage[stretch=10]{microtype}

\usepackage[affil-it]{authblk}

\usepackage{xfrac}

\graphicspath{./img}

\usepackage{hyperref}

\usepackage[nameinlink,noabbrev,capitalize]{cleveref}

\usepackage[plain]{fullpage}

\usepackage{subcaption,graphicx}

\theoremstyle{plain}
\newtheorem{theorem}{Theorem}
\newtheorem*{theorem*}{Theorem}
\newtheorem{proposition}[theorem]{Proposition}

\newtheorem{conjecture}[theorem]{Conjecture}
\newtheorem{lemma}[theorem]{Lemma}
\newtheorem*{lemma*}{Lemma}

\theoremstyle{definition}
\newtheorem{definition}[theorem]{Definition}

\theoremstyle{remark}

\newtheorem{remark}[theorem]{Remark}
\newtheorem*{remark*}{Remark}			

\date{\today}
\title{\textbf{\uppercase{\large{Distortion coefficients \\ of the $\alpha$-Grushin plane}}}}
\author{Samuël Borza\thanks{\href{mailto:samuel.borza2@durham.ac.uk}{samuel.borza2@durham.ac.uk}} }
\affil[]{Department of Mathematical Sciences, Durham University}					

\begin{document}

\maketitle							

\providecommand{\keywords}[1]
{
	\textbf{\textit{Keywords---}} #1
}

\providecommand{\msc}[1]
{
	\textbf{\textit{MSC (2020)---}} #1
}

\begin{abstract}
	 We compute the distortion coefficients of the $\alpha$-Grushin plane. They are expressed in terms of generalised trigonometric functions. Estimates for the distortion coefficients are then obtained and a conjecture of a synthetic curvature condition for the generalised Grushin planes is suggested.
\end{abstract}
\\
\msc{53C17, 51F99, 53B99, 28A75, 28C10}				

\tableofcontents					

\section{Introduction} \label{intro}

Grushin structures first appeared in the work of Grushin on hypoelliptic operators in the seventies, for example, see \cite{grushinoriginal}. The $\alpha$-Grushin plane, denoted by $\mathds{G}_\alpha$, consists of equipping the two-dimensional Euclidean space with the sub-Riemannian structure generated by the global vector fields $X = \partial_x$ and $Y_\alpha = |x|^\alpha \partial_y$.

These structures form a class of rank-varying sub-Riemannian manifolds. In this work, we will focus on the case $\alpha \geqslant 1$. The $\alpha$-Grushin plane has Hausdorff dimension $\alpha + 1$ and is not bracket-generating unless $\alpha$ is an integer. Furthermore, the $\alpha$-Grushin planes constitute a natural generalisation of the traditional Grushin plane, corresponding to the case $\alpha = 1$. Along with the Heisenberg groups $\mathds{H}_n$, they are considered as fundamental examples of sub-Riemannian geometry, exhibiting key characteristics of the theory.

Since the work first set out by Juillet in \cite{heisjuillet} and extended by the same author in \cite{juillet2020subriemanniann}, it is known that, unlike Riemannian manifolds, no sub-Riemannian manifold satisfies the curvature-dimension conditions introduced by Sturm, Lott and Villani. It has been shown by Barilari and Rizzi in \cite{subriemint} that they can, however, support interpolation inequalities and even a geodesic Brunn--Minkowski inequality. For the Heisenberg group, this was in fact first proved by Balogh, Kristály and Sipos in \cite{baloghheis}. Distortion coefficients, which capture some curvature information, play a key role in these results. The present work studies the distortion coefficients of the $\alpha$-Grushin plane.

To achieve this goal, it will be important to study the geodesics of $\mathds{G}_\alpha$ in depth. Because of the lack of a natural connection in sub-Riemannian geometry, geodesics are obtained with Pontryagin's maximum principle. This is the Hamiltonian point of view: a normal minimising path between two points can be lifted to one on the cotangent bundle that satisfies Hamilton's equations. The geodesics of the $\alpha$-Grushin plane were first studied by Li and Chang in \cite{ligrushin}. They are expressed with a generalisation of trigonometric functions, defined as inverses of some special functions. \cref{gentrigunc} and \cref{geodagrushin} are devoted to these topics while in \cref{cutlocusagrushin}, we use an extended Hadamard technique to find the cut loci of $\mathds{G}_\alpha$. The notation $\mathrm{Cut}(q_0)$ stands for the set of cut loci of $q_0$, i.e. the set of points in $\mathds{G}_\alpha$ where the geodesics starting at $q_0$ stop being minimising.

\begin{theorem}[Distortion coefficients of the $\alpha$-Grushin plane] \label{thm:distcoef}
	Let $q_0$ and $q$ be two points of $\mathds{G}_\alpha$ such that $q \notin \mathrm{Cut}(q_0)$. For all $t \in \interval{0}{1}$, we have
	\[
	\beta_t(q_0, q) = \dfrac{\mathrm{J}(t, x_0, u_0, v_0)}{\mathrm{J}(1, x_0, u_0, v_0)},
	\]
	with
	\begin{equation}
	\label{eq:jacdetxuv}
	\mathrm{J}(t, x_0, u_0, v_0) := t \left[u_0 x(t) - (u_0 t + x_0) u(t)\right],
	\end{equation}
	and where $\gamma(t) := (x(t), y(t)) : \interval{0}{1} \to \mathds{G}_\alpha$ denotes the unique constant speed minimising geodesic joining $q_0 = (x_0, y_0)$ to $q$ and $u(t) \mathrm{d}x|_{\gamma(t)} + v(t) \mathrm{d}y|_{\gamma(t)} \in \mathrm{T}_{\gamma(t)}^*(\mathds{G}_\alpha)$ is the corresponding cotangent lift with initial covector $u_0 \mathrm{d}x|_{q_0} + v_0 \mathrm{d}y|_{q_0}$.
\end{theorem}

Because of the analyticity of the geodesic flow, the case $v_0 = 0$ can be seen as taking the limit of $\beta_t(q_0,q)$ as $v_0$ tends to $0$. Geometrically, this means that the points $q_0$ and $q$ are joined by a straight horizontal line. 

\begin{proposition}
	Let $q_0$ and $q$ be two points of $\mathds{G}_\alpha$ such that $q \notin \mathrm{Cut}(q_0)$. When $v_0 = 0$, we have
	\[ 
	\beta_t(q_0,q) = t \dfrac{ (u_0 t + x_0)^{2 \alpha}(u_0 t + x_0) - x_0^{2 \alpha}x_0 }{(u_0 + x_0)^{2 \alpha}(u_0 + x_0) - x_0^{2 \alpha}x_0},
	\]
	for all $t \in \interval{0}{1}$.
\end{proposition}

Although the $\mathrm{CD}$ condition is not suited to this type of spaces, the weaker measure contraction property introduced independently by Ohta and Sturm in \cite{ohtamcp} and \cite{sturm2} seems more adapted to sub-Riemannian geometry. Indeed, there are numerous examples of sub-Riemannian manifolds that do satisfy a $\mathrm{MCP}$ condition, including the Heisenberg group $\mathds{H}_n$ (see \cite{heisjuillet}) and the Grushin plane $\mathds{G}_1$ (see \cite{subriemint}). We, therefore, investigate the $\mathrm{MCP}$ condition for $\mathds{G}_\alpha$ and we obtain a relevant estimate on the distortion coefficients for singular points, that is to say, those on the $y$-axis, and for those lying on the same horizontal line. We therefore propose the following conjecture.

\begin{conjecture}[Curvature dimension of the $\alpha$-Grushin plane]
	\label{conjMCP}
	For $\alpha \geqslant 1$, the $\alpha$-Grushin plane satisfies the measure contraction property condition $\mathrm{MCP}(K, N)$ if and only if $K \leqslant 0$ and
	\[
	N \geqslant 2 \left[ \dfrac{(\alpha + 1) m_\alpha + 1}{m_\alpha + 1} \right]
	\]
	with $m_\alpha \in \interval{-3}{-2}$ the unique non-zero solution of 
	\[
	(m + 1)^{2 \alpha}(m + 1) - \left( (2 \alpha + 1) m + 1 \right) = 0.
	\]
\end{conjecture}

We will provide evidence in favour of this conjecture in \cref{mcpagrushin}. As we will see, the $\mathrm{MCP}(0, N)$ condition is equivalent to a lower bound for the distortion coefficients of the form $\beta_t(q_0, q) \geqslant t^N$. It will be proven that the lower bound holds for singular points. Furthermore, it seems sharp for the points lying on the same horizontal line.

Note that in this work, we always write $(\cdot)^{2 \alpha}$ for $((\cdot)^{2})^\alpha$ and a subscript will sometimes denote a partial derivative.			

\section*{Acknowledgements}

I would like to thank Prof Wilhelm Klingenberg for his supervision and precious encouragements. I am also grateful to Prof Nicolas Juillet, Prof Alpár Mészáros, Prof Alexandru Kristály and Dr Kenshiro Tashiro for reading this work and for providing invaluable feedback. I benefited greatly from instructive discussions with Prof Shingo Takeuchi on the topic of generalised trigonometry. I am thankful to Prof Li Yutian for introducing me to the study of the $\alpha$-Grushin plane and for his constant support. 

This work was supported by the UK Engineering and Physical Sciences Research Council (EPSRC) grant EP/N509462/1 (1888382) for Durham University.		

\section{Preliminaries} \label{prelim}

\subsection{Synthetic curvature-dimension conditions} \label{cdcond}

In this section, we give an overview of metric geometry, synthetic notions of curvature and distortion coefficients. A metric space $(X, \mathrm{d})$ is a length space if the distance is induced from a length structure. This means that $\mathrm{d}(x, y) := \inf \{ \mathrm{L}(\gamma) | \gamma : \interval{a}{b} \to X \text{ is admissible, } \gamma(a) = x \text{ and } \gamma(b) = y \}$ where $\mathrm{L} : \mathcal{A} \to \mathds{R} \cap \{ +\infty \}$ is a length functional on a set of admissible paths $\mathcal{A} \subseteq \mathcal{C}_{\mathrm{d}}(X)$. A minimising geodesic is an admissible path $\gamma : \interval{a}{b} \to X$ in $\mathcal{A}$ such that $\mathrm{d}(\gamma(a), \gamma(b)) = \mathrm{L}(\gamma)$. We refer to \cite{burago2001course} for more on metric geometry.

If the space has the property that every two points can be joined by a minimising geodesic that has constant speed, we will say that $(X, \mathrm{d})$ is a geodesic space. When the metric space $(X, \mathrm{d})$ is equipped with a Radon measure $\mathfrak{m}$, the structure $(X, \mathrm{d}, \mathfrak{m})$ is called a geodesic metric measure space. The notion of distortion coefficients fits into this context.

\begin{definition}
	Let $x, y \in X$. The distortion coefficient from $x$ to $y$ at time $t \in \interval{0}{1}$ is
	\begin{equation}
		\label{distortioncoef}
		\beta_t(x, y) = \limsup_{r \to 0^+} \frac{\mathfrak{m}(\mathrm{Z}_t(x, B_r(y)))}{\mathfrak{m}(B_r(y))},
	\end{equation}
	where $\mathrm{Z}_t(x, B_r(y))$ stands for the set of $t$-intermediate points from $x$ to the ball centred at $y$ of radius $r$;
	\[
	\mathrm{Z}_t(A, B) := \{ \gamma(t) | \gamma \in \mathrm{Geo}(X), \gamma(0) \in A \text{ and } \gamma(1) \in B \}
	\]
	whenever $A$ and $B$ are $\mathfrak{m}$-measurable subsets of $X$.
\end{definition}

Note that $\mathrm{Z}_t(x, B_r(y))$ may not be measurable. If this is the case, the measure $\mathfrak{m}$ in the numerator of \eqref{distortioncoef} is understood as the outer measure of $\mathfrak{m}$.

There is an intuitive physical interpretation of the distortion coefficients (quoted from \cite[Chapter 14.]{villaniold}):
\begin{displayquote}%
	$\left[\beta_t(x, y)\right]$ compares the volume occupied by the light rays emanating from the light source $\left[x\right]$, when they arrive close to $\gamma(t)$, to the volume that they would occupy in a flat space.
\end{displayquote}
In particular, we can thus heuristically expect that the distortion coefficients are related to the curvature of the space.

The theory of synthetic curvature was developed by Lott, Sturm, and Villani (see \cite{villani1}, \cite{sturm1}, and \cite{sturm2}).  Here we summarise some of the points from their works. We denote by $\mathcal{P}(X)$ the set of Borel probability measures and by $\mathcal{P}^2(X)$ the subset of those with finite second moment. We write $\mathrm{Geo}(X)$ for the set of all minimising geodesics of $X$ parametrised by constant speed on $\interval{0}{1}$. For all $t \in \interval{0}{1}$, the evaluation map is defined as
\[
\mathrm{e}_t : \mathrm{Geo}(X) \to X : \gamma \mapsto \gamma(t). 
\]
A dynamical transference plan $\Pi$ is a Borel probability measure on $\mathrm{Geo}(X)$ while a displacement interpolation associated to $\Pi$ is a path $(\mu_t)_{t \in \interval{0}{1}} \subseteq \mathcal{P}^2(X)$ such that $\mu_t = (\mathrm{e}_t)_\# \Pi$ for all $t \in \interval{0}{1}$.
We equip $\mathcal{P}^2(X)$ with the $L_2$-Wasserstein distance $\mathcal{W}_2$: for any $\mu, \nu \in \mathcal{P}^2(X)$, 
\[
\mathcal{W}_2(\mu_0, \mu_1) := \inf_{\pi \in \Pi(\mu, \nu)} \int_X \mathrm{d}(x, y)^2 \pi(\mathrm{d}x\mathrm{d}y),
\]
with $\Pi(\mu, \nu) := \{ \omega \in \mathcal{P}(X^2) | (\mathrm{proj_1})_\# \omega = \mu \text{ and } (\mathrm{proj_2})_\# \omega = \nu \}$. For $\mu_0, \mu_1 \in \mathcal{P}^2(X)$, the set $\mathrm{OptGeo}(\mu_0, \mu_1)$ is the space of all measures $\nu \in \mathcal{P}(\mathrm{Geo}(X))$ such that $(e_0, e_1)_\# \nu$ realises the minimum for the $L_2$-Wasserstein distance. A measure $\nu \in \mathrm{OptGeo}(\mu_0, \mu_1)$ is called a dynamical optimal plan. We now need to define the distortion coefficients of the $(K, N)$-model space. For $K \in \mathds{R}, N \in \interval{1}{+\infty}, \theta \in \interval[open]{0}{+\infty}$ and $t \in \interval{0}{1}$, we set
\[
\tau_{K, N}^{(t)}(\theta) = t^{1/N} \sigma_{K,N - 1}^{(t)}(\theta)^{1 - 1/N},
\]
with
\[
\sigma_{K,N}^{(t)}(\theta) =
\begin{cases}
+\infty & K \theta^2 \geqslant N \pi^2 \\
\dfrac{\sin(t \theta \sqrt{K/N})}{\sin(\theta \sqrt{K/N})} & \text{if } 0 < K \theta^2 < N \pi^2 \\
t & \text{if } K \theta^2 < 0 \text{ and } N = 0 \text{ or if } K \theta^2 = 0 \\
\dfrac{\sinh(t \theta \sqrt{-K/N})}{\sinh(\theta \sqrt{-K/N})} & \text{if } K \theta^2 \leqslant 0 \text{ and } N > 0
\end{cases}.
\]

The definition of the coefficients $\tau_{K, N}$  is not arbitrary. In fact, they are nothing but the distortion coefficients of the model space $X_{(K, N)}$; that is to say, $X_{(K, N)}$ is the $N$-sphere of constant cuvature $K$ if $K > 0$, $X_{(K, N)}$ is the $N$-Euclidean space if $K = 0$ and $X_{(K, N)}$ is the $N$-hyperbolic plane of constant curvature $K$ if $K < 0$ (see again \cite[Chapter 14.]{villaniold}).

We are ready to introduce a first notion of synthetic curvature: the curvature-dimension condition.

\begin{definition}
	Let $K \in \mathds{R}$ and $N \in \interval[open right]{1}{+\infty}$. A geodesic metric measure space $(X, \mathrm{d}, \mathfrak{m})$ satisfies $\mathrm{CD}(K, N)$ if, for any $\mu_0, \mu_1 \in \mathcal{P}^2(X, \mathfrak{m})$ with bounded support, there exists $\nu \in \mathrm{OptGeo}(\mu_0, \mu_1)$ and a $\mathcal{W}_2$-optimal plan $\pi \in \mathcal{P}(X^2)$ such that $\mu_t := (\mathrm{e}_t)_\# \nu \ll \mathfrak{m}$ and for any $N' \geqslant N$,
	\[
	\mathcal{E}_{N'}(\mu_t) \geqslant \int_{X^2} \tau_{K,N}^{(1 - t)}(\mathrm{d}(x, y)) \rho_0^{-1/N'} + \tau_{K,N}^{(t)}(\mathrm{d}(x, y)) \rho_1^{-1/N'} \pi(\mathrm{d}x\mathrm{d}y),
	\]
	where $\mathcal{E}_N$ stands for the Rényi functional
	\[
	\mathcal{E}_N : \mathcal{P}(X) \to \interval{0}{+ \infty} : \rho \mathfrak{m} + \mu_s \mapsto \int_X \rho^{1 - 1/N} \mathfrak{m}(\mathrm{d}x).
	\]
\end{definition}

For an extensive treatment of the $\mathrm{CD}$-condition and more generally of optimal transport theory, we refer the reader to \cite{villaniold}. Alongside this notion of curvature, a weaker condition was developed independently by Sturm and Ohta: the measure contraction property (see \cite{sturm2}, \cite{ohtamcp}).

\begin{definition}
	Let $K \in \mathds{R}$ and $N \in \left[ 1, +\infty \right)$. A geodesic metric measure space $(X, \mathrm{d}, \mathfrak{m})$ satisfies $\mathrm{MCP}(K, N)$ if, for every $x \in X$ and measurable set $A \subseteq X$ with $\mathfrak{m}(A) \in \left( 0, +\infty \right)$, there exists $\nu \in \mathrm{OptGeo}\left( \mu_A, \delta_x \right)$ such that for all $t \in \interval{0}{1}$
	\[
	\mu_A \geqslant (\mathrm{e}_t)_\# \left( \tau_{K, N}^{(1 - t)}(\mathrm{d}(\gamma(0), \gamma(1))) \nu(\mathrm{d}\gamma) \right),
	\]
	where $\mu_A := \frac{1}{\mathfrak{m}(A)} \mathfrak{m} \in \mathcal{P}(X)$ is the normalisation of $\mu|_A$.
\end{definition}

Both the CD and MCP conditions generalise the notion of Ricci curvature
bounded from below by $K \in \mathds{R}$ and dimension bounded from above by $N \geqslant 1$ from Riemannian geometry. Indeed, if $(M, g)$ is a Riemannian manifold and $\psi$ a positive $\mathcal{C}^2$ function on $M$, $\mathrm{d}_g$ the Riemannian distance and $\mathrm{vol}_g$ the Riemannian smooth volume, then the metric measure space $(M, \mathrm{d}_g, \psi \cdot \mathrm{vol}_g)$ satisfies the $\mathrm{CD}(K,N)$ condition if and only if $\dim(M) \leqslant N$ and if $\mathrm{Ric}_{g, \psi, N} \geqslant K g$ where 
\[
\mathrm{Ric}_{g, \psi, N} := \mathrm{Ric_g} - (N - n) \frac{\nabla^2_g \psi^{\frac{1}{N - n}}}{h^{\frac{1}{N - n}}}.
\]
Note that in the case where $N = n$, it only makes sense to consider constant functions $\psi$ in the definition of the generalised Ricci tensor. The proof of the equivalence with the $\mathrm{CD}$ condition can be found in \cite{sturm2} and \cite{villani1}. Furthermore, it is also proved in \cite[Theorem 3.2., Corollary 3.3.]{ohtamcp} that, in the Riemannian setting, the $\mathrm{MCP}(K, N)$ condition is equivalent to $\mathrm{CD}(K, N)$ if $N$ is greater than the topological dimension of $(M, g)$.

For general metric measure spaces, the two notions of synthetic curvature are not equivalent. However, the $\mathrm{CD}$ condition does imply the $\mathrm{MCP}$ condition when the space is \textit{non-branching} (see \cite{sturmcdmcp}). As we will see later, this already appears in sub-Riemannian geometry.

\subsection{Sub-Riemannian geometry} \label{subriemsec}

In what follows, we set up the basics of sub-Riemannian geometry. We rely on \cite{comprehensivesubr} for the general theory (see also \cite{geocontbullo} for vector fields that are not necessary of class $C^\infty$). 

A manifold is a set equipped with an equivalence class of differentiable atlases such that its manifold topology is connected, Hausdorff and second-countable. Here we emphasise the theory of sub-Riemannian manifolds of class $\mathcal{C}^r$ instead of class $\mathcal{C}^{\infty}$. As we will see later, the $\alpha$-Grushin plane is a sub-Riemannian manifold that is generated by global vector fields that might not be smooth.

\begin{definition}
Let $M$ be a smooth manifold of class $\mathcal{C}^r$ for $r \in \mathds{N}^{\geqslant 1} \cup \{ \infty \} \cup \{ \omega \}$. A triple $(E, \langle \cdot, \cdot \rangle_E, f_E)$ is said to be a sub-Riemannian structure of class $\mathcal{C}^r$ on $M$ if
\begin{enumerate}
\item $E$ is a $\mathcal{C}^r$-vector bundle on $M$,
\item $\langle \cdot, \cdot \rangle_E$ is a $\mathcal{C}^r$-Euclidean metric on $E$,
\item $f_E : E \to \mathrm{T}(M)$ is a $\mathcal{C}^r$-morphism of vector bundles.
\end{enumerate}
\end{definition}
The family $\mathcal{D}$ of $\mathcal{C}^r$-horizontal vector fields is defined as
\[
\mathcal{D} := \{ f_E \circ u | u \text{ is a section of } E \text{ of class } \mathcal{C}^r \}.
\]
We also define the distribution at point a $p \in M$ with
\[
\mathcal{D}_p := \{ v(p) \mid v \in \mathcal{D} \}.
\]
The rank of the sub-Riemannian structure at $p \in M$ is $\mathrm{rank}(p) := \mathrm{dim}(\mathcal{D}_p)$. Observe that in our definition, a sub-Riemannian manifold can be rank-varying; i.e. the map $\mathrm{rank}(\cdot)$ might not be constant.

\begin{definition}
	We say that curve $\gamma : \interval{0}{T} \to M$ is horizontal if $\gamma$ is Lipschitz in charts and if there exists a control $u \in \mathrm{L}^2(\interval{0}{T}, E)$ such that for all $t \in \interval{0}{T}$, we have $u(t) \in E_{\gamma(t)}$ and $\dot{\gamma} (t) = f_E(u(t))$.
	The sub-Riemannian length of $\gamma$ is defined by
	\[
	\mathrm{L}_{\mathrm{CC}}(\gamma) = \int_{0}^{T} \lVert \dot{\gamma} (t) \rVert_{\gamma(t)} \mathrm{d}t,
	\]
	where $\lVert v \rVert_{\mathcal{D}_p} := \min \left\{ \sqrt{g(u, u)} \mid u \in E_p \text{ and } f_E(u) = (p, v) \right\}$ for $v \in \mathcal{D}_p$ and $p \in M$.
\end{definition}

\begin{remark}
	It can be proven that $\lVert \cdot \rVert_{\mathcal{D}_p}$ is well defined, induced by an inner product $\langle \cdot, \cdot \rangle_{\mathcal{D}_p}$ and that the map $\lVert \dot{\gamma} (\cdot) \rVert_{\mathcal{D}_{\gamma(\cdot)}}$ is measurable.
\end{remark}

In the case where every two points can be joined by a horizontal curve, we have a well-defined distance function on $M$.

\begin{definition}
	Let $M$ be a sub-Riemannian manifold. The sub-Riemannian distance $\mathrm{d_{CC}}$ of $M$, also called the Carnot--Carathéodory distance, is defined by
	\[
	\mathrm{d_{CC}}(x, y) := \inf \{ \mathrm{L_{CC}}(\gamma) | \gamma : [0,T] \to M \text{ is horizontal and } \gamma(0) = x \text{ and } \gamma(T) = y \}.
	\]
\end{definition}

Traditionally, the definition of a sub-Riemannian structure demands that $\mathcal{D}$ is a $\mathcal{C}^\infty$-distribution and that it satisfies the Hörmander condition; that is to say, $\mathrm{Lie}_p(\mathcal{D}) = \mathrm{T}_p(M)$ for all $p \in M$. This is motivated by the following well-known result.

\begin{theorem}[Chow-–Rashevskii theorem]
	Let $M$ be a sub-Riemannian manifold such that its distribution $\mathcal{D}$ is $\mathcal{C}^\infty$ and satisfies the Hörmander condition. Then, $(M, \mathrm{d_{CC}})$ is a metric space and the manifold and metric topology of $M$ coincide.
\end{theorem}

We refrain from this convention here, as the Grushin planes that we will study do not always satisfy this property. However, we will assume hereafter that every two points of the sub-Riemannian manifold $M$ can be joined by a horizontal curve, making $\mathrm{d_{CC}}$ a distance of $M$, and that the metric and manifold topologies do coincide.

Finally, the horizontal distribution of a sub-Riemannian manifold $M$ is defined by 
\[
\mathrm{H}(M) := \bigsqcup_{p \in M} \mathcal{D}_p .
\]
Note that $\mathrm{H}(M)$ has no natural structure of subbundle in $\mathrm{T}(M)$ if $M$ is rank-varying.

Now that we have turned our sub-Riemannian manifold into a metric space, we would like to study the geodesics associated with $\mathrm{d_{CC}}$. These would be horizontal curves that are locally a minimiser for the length functional $\mathrm{L_{CC}}$. Because of the lack of a torsion-free metric connection, we cannot study geodesics through a covariant derivative. Rather, some sub-Riemannian geodesics can be characterised via Hamilton's equation.

Given $m$-global $\mathcal{C}^r$-vector fields $X_1, \dots, X_m : M \to \mathrm{T}(M)$ on a $\mathcal{C}^r$-manifold $M$, we can induce on $M$ a sub-Riemannian structure in the following way. We set $E = M \times \mathds{R}^m$ the trivial bundle of rank $m$, $f_E : E \to \mathrm{T}(M) : (p, (u_1, \dots, u_m)) \mapsto \sum_{k = 1}^{m} u_k X_k(p)$ and finally the metric on $E$ is the Euclidean one. In this way, we induce an inner product on $\mathcal{D}_p = \mathrm{span}\{X_1(p), ..., X_m(p)\}$ by the polarisation formula applied to the norm
\begin{equation}
	\label{innerproduct}
	\|u\|^2_{\mathcal{D}_p} := \min \left\{ \sum_{k = 1}^{m} u^2_i \mid \sum_{k = 1}^{m} u_i X_k(p) = u \right\}.
\end{equation}
The family $(X_1, \dots, X_m)$ is said to be a generating family of the sub-Riemannian manifold $(M, E, \langle \cdot, \cdot \rangle_E, f_E)$. A \textit{free} sub-Riemannian structure is one that is induced from a generating family. Every sub-Riemannian structure is \textit{equivalent} to a free one (see \cite[Section 3.1.4]{comprehensivesubr}). From now on, we will therefore assume that the sub-Riemannian manifolds considered are free.

\begin{definition}
	Let $M$ be a sub-Riemannian manifold and $(X_1, \dotsc, X_m)$ a generating family of vector fields. The Hamiltonian of the sub-Riemannian structure is defined by
	\[
	H : \mathrm{T^*}(M) \to \mathds{R} : (p, \lambda_0) \mapsto H(p, \lambda_0) := \frac{1}{2}\sum_{k = 1}^m h_k(p, \lambda_0)^2,
	\]
	where $h_k(p, \lambda_0) := \langle \lambda_0, X_k(p) \rangle$.
\end{definition}

We therefore approach the problem via the cotangent bundle $\mathrm{T}^*(M)$, on which there is a natural symplectic form $\sigma$. We can now characterise length minimisers of a sub-Riemannian manifold.

\begin{theorem}[Pontryagin's maximum principle]
	Let $\gamma : \interval{0}{T} \to M$ be a horizontal curve which is a length minimiser parametrised by constant speed. Then, there exists a Lipschitz curve $\lambda(t) \in \mathrm{T}^*_{\gamma(t)}(M)$ such that one and only one of the following is satisfied:
	\begin{enumerate}
		\item[(N)] $\dot{\lambda} = \overrightarrow{H}(\lambda)$, where $\overrightarrow{H}$ is the unique vector field in $\mathrm{T}^*(M)$ such that $\sigma(\cdot, \overrightarrow{H}(\lambda)) = \mathrm{d}_\lambda H$ for all $\lambda \in \mathrm{T}^*(M)$;
		\item[(A)] $\sigma_{\lambda(t)}(\dot{\lambda} (t), \cap_{k = 1}^n \mathrm{ker}(\mathrm{d}_{\lambda(t)} h_k)) = 0$ for all $t \in \interval{0}{T}$.
	\end{enumerate}
\end{theorem}

If $\lambda$ satisfies $(N)$ (resp. $(A)$), we will say that $\lambda$ is a normal extremal (resp. abnormal extremal) and $\gamma$ is a normal geodesic (resp. abnormal geodesic). Note that a geodesic may be both normal and abnormal. The projection of a normal extremal onto $M$ is locally minimising, that is to say a (normal) geodesic parametrised by constant speed.

If $\gamma$ is a normal geodesic associated with a normal extremal $\lambda$, then $(N)$ is nothing but Hamilton's equation for $H$ in the natural coordinates of the cotangent bundle:
\begin{equation}
\label{eq:hamiltonode}
\left\{
\begin{array}{rcl}
\dot{x}_i & = & \dfrac{\partial H}{\partial p_i} \\ 
\dot{p}_i & = & -\dfrac{\partial H}{\partial x_i}. \\
\end{array}
\right.
\end{equation}

The \textit{exponential map} at $p \in M$ is the function
\[
\mathrm{exp}_p : \mathscr{A}_p \to M : \lambda_0 \mapsto \pi(\mathrm{e}^{\vec{H}}(\lambda_0)),
\]
where $\pi : \mathrm{T}^*(M) \to M$ is the projection, $\mathrm{e}^{t\overrightarrow{H}}$ is the flow of $\overrightarrow{H}$ and $\mathscr{A}_p \subseteq \mathrm{T}^*_p(M)$ is the open set of covectors such that the corresponding solution of \eqref{eq:hamiltonode} is defined on the whole interval $\interval{0}{1}$.

The cut time of a geodesic $\gamma$ is defined as
\[
t_{\mathrm{cut}}[\gamma] := \sup \{ t > 0 \mid \gamma|_{\interval{0}{t}} \text{ is minimising} \}.
\]
When $t_{\mathrm{cut}}[\gamma] < +\infty$, we say that $\gamma(t_{\mathrm{cut}}[\gamma])$ is the cut point to $\gamma(0)$ along $\gamma$. If $t_{\mathrm{cut}}[\gamma] = +\infty$, we say that $\gamma$
has no cut point. We denote by $\mathrm{Cut(q_0)}$ the set of all cut points of geodesics starting from a point $q_0 \in M$.

The study of abnormal geodesics is an area of intensive research. It does happen that a sub-Riemannian structure does not have any non-trivial abnormal geodesic (the trivial geodesic is always abnormal as soon as $\mathrm{rank}(\mathcal{D}_p) < \mathrm{dim} (M)$). In this case, a sub-Riemannian manifold is said to be \textit{ideal}.

The $\mathrm{CD}(K, N)$ condition is never satisfied for ideal sub-Riemannian manifolds $M$ such that $\mathrm{rank}(\mathcal{D}_p) < \mathrm{dim} (M)$ at every point $p \in M$ (see \cite{juillet2020subriemanniann}). However, it is known that they often satisfy an $\mathrm{MCP}$ condition: the Heisenberg groups (see \cite{heisjuillet}), generalised H-type groups, Sasakian manifolds (see \cite[Section 7.]{subriemint}), etc. We conclude this section with the following theorem that relates the $\mathrm{MCP}$ condition to a lower bound on the distortion coefficients of an ideal sub-Riemannian manifold.

\begin{theorem}[{\cite[Theorem 9.]{subriemint}}]
	\label{equivMCP}
	Let $M$ be an ideal sub-Riemannian manifold equipped with a smooth measure $\mu$. When $N \geqslant 1$, the following conditions are equivalent:
	\begin{enumerate}[label=\normalfont(\roman*)]
		\item $\beta_t(q_0, q) \geqslant t^N$ for all $q_0, q \notin \mathrm{Cut}(M)$ and $t \in [0, 1]$;
		\item The measure contraction property $\mathrm{MCP}(0, N)$ is satisfied, i.e. for all non-empty Borel sets $B \subseteq M$ and $q \in M$ we have $\mu(\mathrm{Z}_t(q, B)) \geqslant t^N \mu(B)$. 
	\end{enumerate}
\end{theorem}

\subsection{Generalised trigonometric functions} \label{gentrigunc}

In this section, we give an account of $(p, q)$-trigonometry. The generalised sine and cosine functions will be essential in the study of the geometry of the $\alpha$-Grushin plane, as shown by Li in \cite{ligrushin}. Generalised trigonometry has a long history. The theory as presented here was pioneered by Edmunds in \cite{edmunds2}. For recent developments, we point out the work of Takeuchi \cite{takeuchi1} and the references therein, as well as \cite{convexlok} for a related approach via convex geometry.

Consider
\[
F_{p, q} : \interval{0}{1} \to \mathds{R} : x \mapsto \int_0^x \dfrac{1}{\sqrt[p]{1 - t^{q}}} \mathrm{d}t.
\]
The map $F_{p,q}$ being strictly increasing, we may define its inverse
\[
\sin_{p,q} : \interval[scaled]{0}{\dfrac{\pi_{p,q}}{2}} \to \mathds{R} : x \mapsto F_{p,q}^{-1}(x),
\]
where the $(p, q)$-pi constant is defined as
\[ 
\pi_{p,q} := 2 \int_0^1 \dfrac{1}{\sqrt[p]{1 - t^{q}}} \mathrm{d}t = \mathrm{B}\left( \dfrac{1}{p}, 1 - \dfrac{1}{q} \right).
\]
Here the function $\mathrm{B}(\cdot, \cdot)$ stands for the complete beta function.

We will extend the $(p, q)$-sine function to the whole real line. We first note that $\sin_{p,q}(0) = 0$ and $\sin_{p,q}(\pi_{p,q}/2)$ = 1. For $x \in \interval[scaled]{\pi_{p,q}/2}{\pi_{p,q}}$, we set $\sin_{p,q} (x) := \sin_{p,q}(\pi_{p,q} - x)$. The $(p, q)$-sine is then extended to $\interval{-\pi_{p, q}}{\pi_{p, q}}$ by requiring that it is odd and finally to the whole $\mathds{R}$ by $2\pi_{p,q}$-periodicity. We then define the $(p, q)$-cosine by setting $\cos_{p,q} := (\sin_{p,q})'$. These two functions are of class $\mathcal{C}^1$. In fact, they are also of class $\mathcal{C}^{\infty}$ except at the points $x = k \pi_{p,q}$ for $k \in \mathds{Z}$.

We have the following identities:
\begin{equation}
	\left\{ 
	\begin{array}{l}
		|\sin_{p,q}|^q + |\cos_{p,q}|^p = 1, \\
		(\sin_{p,q})'' =  (\cos_{p,q})' = \dfrac{-q}{p} |\cos_{p,q}|^{2 - p} |\sin_{p,q}|^{q - 2} \sin_{p,q}.
	\end{array}
	\right.
\end{equation}
Therefore, the $(p, q)$-sine function can be alternatively defined as the solution to the following ordinary differential equation
\begin{equation}
	\label{eqn:odegrushin}
	-(|f'|^{p - 2} f')' = \dfrac{(p - 1)q}{p} |f|^{q - 2} f, \ \ f(0) = 0, \ \ f' (0) = 1.
\end{equation}

As for the usual sine and cosine functions, we have $\sin_{p,q}(x + \pi_{p,q}) = - \sin_{p,q}(x)$ and $\cos_{p,q}(x + \pi_{p,q}) = -\cos_{p,q}(x)$. However, unlike the case of classical trigonometric functions, general addition formulas are not known for $\sin_{p,q}(x + y)$ and $\cos_{p,q} (x + y)$ (except for very specific values of $p$ and $q$). This problem ultimately comes down to finding a function $F_{p, q}$ that solves the integral equation
\[
\int_0^{F_{p, q}(x, y)} \dfrac{1}{\sqrt[p]{1 - t^{q}}} \mathrm{d}t = \int_0^{x} \dfrac{1}{\sqrt[p]{1 - t^{q}}} \mathrm{d}t + \int_0^{y} \dfrac{1}{\sqrt[p]{1 - t^{q}}} \mathrm{d}t.
\]
We would then have $\sin_{p,q}(x + y) = F_{p,q}(\sin_{p,q}(x), \sin_{p,q}(y))$. This is a very difficult problem, even for integer values of $p$ and $q$. For $(p, q) = (2, 2)$, the classical addition formula for the sine functions emerges. When $(p, q) = (2, 4)$, the corresponding addition formula is the one used for the lemniscate function that Euler investigated in \cite{eulerlem}: let $\mathrm{sl}(x) := \sin_{2,4}(x)$ (resp. $\mathrm{sl}'(x) := \cos_{2,4}(x)$) stand for the sinlem function (resp. the sinlem' function), then we have
\[
\mathrm{sl}(x + y) = \dfrac{\mathrm{sl}(x)\mathrm{sl}'(y) + \mathrm{sl}(y)\mathrm{sl}'(x)}{1 + \mathrm{sl}^2(x) \mathrm{sl}^2(y)},
\]
with an analogous formula for $\mathrm{sl}'(x + y)$. Note that Euler's coslem function is defined as $\mathrm{cl(x)} = \mathrm{sl}(x + \pi_{(2, 4)}/2)$, which is different from our $(2, 4)$-cosine function.				

\section{Geometry of the \texorpdfstring{$\alpha$}{alpha}-Grushin plane} \label{geom}

\subsection{Geodesics of the \texorpdfstring{$\alpha$}{alpha}-Grushin plane} \label{geodagrushin}

For $\alpha \in [1, +\infty)$, the $\alpha$-Grushin plane $\mathds{G}_\alpha$ is defined as the sub-Riemannian structure on $\mathds{R}^2$ generated by the global vector fields $X = \partial_x$ and $Y_\alpha = |x|^{\alpha} \partial_y$, as explained in \cref{subriemsec}. This generating family of vector fields is $\mathcal{C}^{\lfloor \alpha \rfloor}$ if $\alpha$ is not an integer and $\mathcal{C}^\infty$ otherwise. 

The horizontal space at $p \in \mathds{G}_\alpha$ is $ \mathcal{D}_p \left( \mathds{G}_\alpha \right) = \mathrm{span} \{ X(p), Y_\alpha(p) \}$ and the horizontal distribution is the disjoint union of these $\mathrm{H}(\mathds{G}_\alpha) = \sqcup_{p \in \mathds{G}_\alpha} H_p(\mathds{G}_\alpha)$. The rank of $\mathcal{D} = \mathrm{span} \left\lbrace X, Y_\alpha \right\rbrace $ is not constant: it is a singular distribution if $x = 0$ and Riemannian otherwise. We then consider the scalar metric $\langle \cdot, \cdot \rangle_{\mathcal{D}_p}$ on $\mathcal{D}_p$ as described in \eqref{subriemsec}. If for example $u X(x, y) + v Y_\alpha(x, y) \in \mathcal{D}_{(x, y)}$ and $x \neq 0$, then
\[
\langle u, v \rangle_{\mathcal{D}_{(x, y)}} = u^2 + \frac{1}{x^{2 \alpha}} v^2.
\]
This turns the $\alpha$-Grushin plane $\mathds{G}_\alpha$ into a sub-Riemannian manifold. It is easy to see that it does not satisfy the Hörmander condition unless $\alpha$ is an integer.

Let $I$ be a non-empty interval of $\mathds{R}$. As we have seen in the previous section, a path $\gamma : I \to \mathds{G}_\alpha$ is said to be horizontal if, for almost every $t \in I$, the equality $\dot{\gamma}(t) = u(t) X(\gamma(t)) + v(t) Y_\alpha(\gamma(t))$ holds for some $L^2$-maps $u, v : I \to \mathds{R}$. In particular, this implies that $\dot{\gamma} (t) \in \mathcal{D}_{\gamma(t)}$ for almost every $t \in I$. We can compute the length of a horizontal curve with the formula $\mathrm{L}_\alpha(\gamma) = \int_I \| \dot{\gamma}(t)\|_{\mathcal{D}_{\gamma(t)}} \mathrm{d}t$. We denote the Carnot--Carathéodory distance associated with $\mathrm{L}_\alpha$ by $\mathrm{d}_\alpha$. Equipping the $\alpha$-Grushin plane with the Lebesgue measure $\mathcal{L}^2$, we obtain a metric measure space $(\mathds{G}_\alpha, \mathrm{d}_\alpha, \mathcal{L}^2)$.

The theory of sub-Riemannian geometry informs us that the geodesics of the space are found by solving Hamilton's equations. Here, the Hamiltonian is
\[
H : \mathrm{T}^*(\mathds{G}_\alpha) \to \mathds{R} : (x, y, u \mathrm{d}x|_{(x, y)} + v \mathrm{d} y|_{(x, y)}) \mapsto \frac{1}{2} (u^2 + v^2 x^{2 \alpha}).
\]
A simple calculation shows that there are no non-trivial abnormal geodesics in the $\alpha$-Grushin plane. Consequently, the sub-Riemannian manifold $\mathds{G}_\alpha$ is ideal. In this context, Hamilton's equations \eqref{eq:hamiltonode} become
\begin{equation}
	\label{eq:odegrushin}
	\left\{
	\begin{array}{l}
		\dot{x} = u, \\
		\dot{y} = v x^{2 \alpha} \\
		\dot{u} = - \alpha v^2 x^{2 (\alpha - 1)}x \\
		\dot{v} = 0
	\end{array}.
	\right.
\end{equation}

We observe that $\ddot{x} = - \alpha v^2 x^{2(\alpha - 1)}$. When $v_0 = 1$, this is just the equation \eqref{eqn:odegrushin} for $(p, q) = (2, 2\alpha)$. The $(2, 2 \alpha)$-trigonometric functions will therefore be essential and in what follows, we will denote $\sin_\alpha$ instead of $\sin_{2,2 \alpha}$ (and respectively $\cos_\alpha$, $\pi_\alpha$) for simplicity.

\begin{theorem}
	\label{thm:geodesics}
	Let $\gamma : I \to \mathds{G}_\alpha$ be a horizontal path with initial value $\gamma(0) = (x_0, y_0)$ and $\lambda(t) = u(t) \mathrm{d}x|_{\gamma(t)} + v(t) \mathrm{d}y|_{\gamma(t)}$ be the cotangent lift with initial covector $(u(0), v(0)) = (u_0, v_0)$. 
	
	In the case where $v_0 \neq 0$ and $(x_0, u_0) \neq 0$, the curve $\gamma$ is a geodesic if and only if 
	\begin{equation}
		\left\{
		\begin{array}{rcl}
			x(t) & = & A \sin_\alpha(\omega t + \phi) \label{geodesics}  \\
			y(t) & = & y_0 + v_0 \dfrac{A^{2 \alpha}}{(\alpha + 1)\omega^2} \Big[ \omega^2 t + \omega \cos_\alpha(\phi) \sin_\alpha(\phi) \\
			& & \qquad \qquad \qquad \qquad \qquad - \omega \cos_\alpha(\omega t + \phi) \sin_\alpha(\omega t + \phi) \Big]\\
			u(t) & = & A \omega \cos_\alpha(\omega t + \phi) \\
			v(t) & = & v_0 \\
		\end{array}
		\right.
	\end{equation}
	for uniquely determined parameters $A, \omega \in \mathds{R} \setminus \{0\}$ and $\phi \in \left[0, 2 \pi_\alpha\right)$ satisfying
	\begin{equation}
		\begin{gathered}
			A \omega > 0, \ \ A^2 \omega^2 = u_0^2 + v_0^2 x_0^{2 \alpha}, \ \ \omega^2 = v_0^2 A^{2(\alpha - 1)}, \\ 
			x_0 = A \sin_\alpha(\phi) \text{ and } u_0 = A \omega \cos_\alpha(\phi).
		\end{gathered}\label{basicrel}
	\end{equation}
	If $v_0 = 0$ or $(x_0, u_0) = 0$, the geodesic is $(x(t), y(t)) = (u_0 t + x_0, y_0)$ with its lift being constant: $(u(t), v(t)) = (u_0, v_0)$.
\end{theorem}

\begin{remark}
	Since the right-hand side of the equation is continuous with respect to the initial condition $v_0$, the normal extremals corresponding to $v_0$ can be obtained by letting $v_0$ tend to $0$ in \eqref{geodesics}.
\end{remark}

\begin{proof}
	The case when $v_0 = 0$ or $(x_0, u_0) = 0$ is straightforward. We assume that $v_0 \neq 0$ and $(x_0, u_0) \neq 0$. For $A, \omega \in \mathds{R} \setminus \{0\}$ such that $A w > 0$ and $\phi \in \interval[open right]{0}{2 \pi_\alpha}$, we have
	\begin{align*}
		(A \sin_\alpha(\omega t + \phi))'' &= (A \omega \cos_\alpha(\omega t + \phi))' \\
		&= - \alpha A \omega^2 \sin_\alpha(\omega t + \phi)^{2(\alpha - 1)}\sin_\alpha(\omega t + \phi) \\
		&= - \alpha \frac{\omega^2}{A^{2(\alpha - 1)}} (A\sin_\alpha(\omega t + \phi))^{2(\alpha - 1)} (A\sin_\alpha(\omega t + \phi)).
	\end{align*}
	By the uniqueness of solutions to the differential equation \eqref{eq:odegrushin}, we get
	\begin{equation}
		\left\{
		\begin{array}{rcl}
			x(t) & = & A \sin_\alpha(\omega t + \phi), \\
			u(t) & = & A \omega \cos_\alpha(\omega t + \phi),
		\end{array}
		\right.
	\end{equation}
	where we set $\omega^2 = v_0^2 A^{2(\alpha - 1)}$, $x_0 = A \sin_\alpha(\phi)$ and $u_0 = A \omega \cos_\alpha(\phi)$. Considering the constant of motion $u^2 + v^2 x^{2 \alpha}$ at $t = 0$ yields
	\[
	u_0^2 + v_0^2 x_0^{2 \alpha} = (A \omega \cos_\alpha(\phi))^2 + \frac{\omega^2}{A^{2(\alpha - 1)}}(A \sin_\alpha(\phi))^{2 \alpha} = A^2 \omega^2.
	\]
	Since $\ddot{x} = - \alpha v_0^2 x^{2(\alpha - 1)}x$, we deduce that $x^{2\alpha} = - x \ddot{x} / \alpha v_0^2$ and thus, integrating by part, we have
	\[
	\int_0^t x^{2 \alpha} = \int_0^t \frac{- x \ddot{x}}{\alpha v_0^2} = \frac{- 1}{\alpha v_0^2} \left( [x \dot{x}]_0^t - \int_0^t (\dot{x})^2 \right) = \frac{- 1}{\alpha v_0^2} \left( [x u]_0^t - \int_0^t u^2 \right).
	\]
	We use the identity $u^2 = A^2 \omega^2 - v^2 x^{2 \alpha}$ to find
	\begin{align*}
		\int_0^t x^{2 \alpha} & = \dfrac{- 1}{\alpha v_0^2} \left(x(t)u(t) - x(0)u(0) - \int_0^t A^2 \omega^2 + \int_0^t v^2 x^{2 \alpha} \right) \\ 
		& = \dfrac{A^2}{\alpha v_0^2} \left(\omega^2 t + \omega \cos_\alpha(\phi) \sin_\alpha(\phi) \right. \\
		& \qquad \left. - \omega \cos_\alpha(\omega t + \phi) \sin_\alpha(\omega t + \phi) - \dfrac{v_0^2}{A^2} \int_0^t x^{2 \alpha} \right).
	\end{align*}
	Finally, we isolate $\int_0^t x^{2 \alpha}$ and  integrate $\dot{y} = v_0 x^{2 \alpha}$ to get
	\begin{align*}
		y(t) & = y_0 + v_0 \dfrac{A^{2 \alpha}}{(\alpha + 1)\omega^2}\Big( \omega^2 t + \omega \cos_\alpha(\phi) \sin_\alpha(\phi) \\
		& \qquad \qquad \qquad \qquad \qquad \qquad - \omega \cos_\alpha(\omega t + \phi) \sin_\alpha(\omega t + \phi) \Big).
	\end{align*}
	
	It remains to prove that there is a one-to-one and continuous correspondence between the variables $(A, \omega, \phi)$ and $(x_0,u_0, v_0)$ via \eqref{basicrel}. Going from $(A, \omega, \phi)$ to $(x_0,u_0, v_0)$ is clear. The other direction is given by
	\begin{equation}
		\label{corresp}
		\left\{
		\arraycolsep=1.4pt\def\arraystretch{1.7}
		\begin{array}{rcl}
			A & = & \displaystyle \mathrm{sgn}(v_0) \left( \frac{u_0^2 + v_0^2 x_0^{2 \alpha}}{v_0^2} \right)^{1/2\alpha}, \\
			\displaystyle \omega & = & \displaystyle v_0 \left( \frac{u_0^2 + v_0^2 x_0^{2 \alpha}}{v_0^2} \right)^{(\alpha - 1)/2\alpha}, \\
			\displaystyle \sin_\alpha(\phi) & = & \mathrm{sgn}(v_0) x_0 \left( \dfrac{v_0^2}{u_0^2 + v^2_0 x_0^{2 \alpha}} \right)^{1/2\alpha}, \\
			\displaystyle \cos_\alpha(\phi) & = & \dfrac{u_0}{(u_0^2 + v^2_0 x_0^{2 \alpha})^{1/2}}.
		\end{array}
		\right.
	\end{equation}
\end{proof}
By differentiating the relations \eqref{basicrel} with respect to $x_0, u_0$ and $v_0$, we find the following useful identities:
\begin{equation}
	\begin{gathered}
		A_{x_0} = \dfrac{1 - \cos_\alpha^2(\phi)}{\sin_\alpha(\phi)}, \quad A_{u_0} = \dfrac{\cos_\alpha(\phi)}{\alpha \omega}, \quad A_{v_0} = \dfrac{-\cos^2_\alpha(\phi) A}{\alpha v_0}; \\
		\phi_{x_0} = \dfrac{\cos_\alpha(\phi)}{A}, \quad \phi_{u_0} = \dfrac{-\sin_\alpha(\phi)}{\alpha \omega A}, \quad \phi_{v_0} = \dfrac{\sin_\alpha(\phi)\cos_\alpha(\phi)}{\alpha v_0}; \\
		\omega_{x_0} = (\alpha - 1) \left( \dfrac{\omega}{A} \right) \left( \dfrac{1 - \cos_\alpha^2(\phi)}{\sin_\alpha(\phi)} \right), \quad \omega_{u_0} = \left(\dfrac{\alpha - 1}{\alpha}\right) \dfrac{\cos_\alpha(\phi)}{A},  \\
		\omega_{v_0} = \dfrac{\omega}{v_0}\left(1 - \left(\dfrac{\alpha - 1}{\alpha}\right) \cos^2_\alpha(\phi) \right).
	\end{gathered}\label{der}
\end{equation}

We mention here the work of Li and Chang (see \cite{ligrushin}). They obtained the geodesics joining every two points in the $\alpha$-Grushin plane by solving the boundary value problem corresponding to the differential equation in \cref{thm:geodesics}. We note that their results are stated for $\alpha \in \mathds{N} \setminus \{0\}$. However, if we carefully define sub-Riemannian manifolds of class $\mathcal{C}^k$ as it was done in \cref{subriemsec} and in \cref{geodagrushin}, we can see that their conclusions remain valid in the case $\alpha \geqslant 1$. In particular, their detailed study of the geodesics was used to derive an expression for the Carnot--Carathéodory distance of $\mathds{G}_\alpha$ between every two points. 

\subsection{Cut locus of the \texorpdfstring{$\alpha$}{alpha}-Grushin plane} \label{cutlocusagrushin}

When we look at the the geodesics of $\mathds{G}_\alpha$, we observe three types of behaviours: the straight horizontal lines corresponding to an initial covector with $v_0 = 0$; the geodesics for which $x_0 = 0$ (called \textit{singular} or \textit{Grushin points}); and those for which $x_0 \neq 0$ (called \textit{Riemannian points}). In this section, we investigate the \textit{sub-Riemannian cut loci} and \textit{times} of the $\alpha$-Grushin plane. The techniques used here were developed in \cite[Section 3.2]{agrachevcut}, \cite[Appendix A]{rizzicut} and \cite[Section 13.5]{comprehensivesubr}.

The case when $v_0 = 0$ is trivial: the corresponding geodesic is a straight horizontal line and is length-minimising for all times. Its cut locus is empty and its cut time is infinite.

We now look at a geodesic $\gamma$ starting from a singular point $x_0 = 0$. Since $A^2 \omega^2 = u_0^2 + v_0^2 x_0^{2 \alpha} = u_0^2 = \kappa^2$, where the positive parameter $\kappa > 0$ is the constant speed of the geodesic $\gamma$, we can parametrise $u_0, v_0$ and the corresponding parameters $A$ and $\omega$ with respect to $t \in \mathds{R}$ and $\beta \in \mathds{R} \setminus \{0\}$:
\[
\begin{gathered}
	u_0 = \pm \kappa, \ \ v_0 = \beta, \ \ \phi = 0 \text{ or } \pi_\alpha, \\
	A = \mathrm{sgn}(\beta)\left(\dfrac{\kappa}{|\beta|}\right)^{1/\alpha} \text{and } \omega = \beta \left(\dfrac{\kappa}{|\beta|}\right)^{\frac{\alpha - 1}{\alpha}} \\
\end{gathered}
\]
The geodesic starting at $(0, y_0)$ can then be written as follows:
\begin{equation}
	\label{geodesicsagrushin0}
	\left\{
	\begin{array}{rcl}
		x^{\pm}(t, \beta) & = & \pm \mathrm{sgn}(\beta)\left(\dfrac{\kappa}{|\beta|}\right)^{1/\alpha} \sin_\alpha \left(\beta \left(\dfrac{\kappa}{|\beta|}\right)^{\frac{\alpha - 1}{\alpha}} t \right) \\
		y(t, \beta) & = &  y_0 + \dfrac{1}{(\alpha + 1)} \left(\dfrac{\kappa}{|\beta|}\right)^{\frac{\alpha + 1}{\alpha}} \left[ \beta \left(\dfrac{\kappa}{|\beta|}\right)^{\frac{\alpha - 1}{\alpha}} t \right. \\
		& & \left. - \cos_\alpha\left(\beta \left(\dfrac{\kappa}{|\beta|}\right)^{\frac{\alpha - 1}{\alpha}} t \right) \sin_\alpha\left(\beta \left(\dfrac{\kappa}{|\beta|}\right)^{\frac{\alpha - 1}{\alpha}} t \right) \right]
	\end{array}
	\right.
\end{equation}
and in the case $\beta = 0$, the system can be interpreted as $x^{\pm}(t, \beta) = \pm \kappa t$ and $y(t, \beta) = y_0$. From \eqref{geodesicsagrushin0}, we see that the geodesic $(x^{+}(\cdot, \beta), y(\cdot, \beta))$ is a reflection of $(x^{-}(t, \beta), y(t, \beta))$ with respect to the $y$-axis. Furthermore, these two intersect at the $y$-axis for the first time when $t = \pi_\alpha / |\omega|$. Therefore, a geodesic $\gamma$ starting at a singular point $(0, y_0)$ must lose its optimality after $t = \pi_\alpha / |\omega|$. The following lemma guarantees the optimality of $\gamma$ when $t \leqslant \pi_\alpha / |\omega|$. It is analogous to the case $\alpha = 1$ (see \cite[Section 13.5.2]{comprehensivesubr}).

\begin{lemma}
	A geodesic $\gamma$ starting at a singular point $(0, y_0) \in \mathds{G}_\alpha$ is minimising when $t \leqslant \pi_\alpha / |\omega|$.
\end{lemma}

\begin{proof}
	Let $(x_1, y_1) := \gamma(t^*)$ for a fixed $t^* \in \interval[open right]{0}{\pi_\alpha / |\omega|}$. From \cite[Theorem 12]{ligrushin}, we know that there is a finite number of geodesics joining the singular point $(0, y_0)$ to a point $(x_1, y_1)$, only one among them being minimising. We claim that there is a unique $\beta \in \mathds{R}$ and unique $t \in \interval[open right]{0}{\pi_\alpha / |\omega|}$ such that $(x^{\pm}(t, \beta), y(t, \beta)) = (x_1, y_1)$. By the symmetries of the $\alpha$-Grushin and since $x_1 = 0$ corresponds to $\gamma$ being a horizontal line, we can assume that $x_1 > 0$ and $y_1 \geqslant y_0$ without loss of generality. In particular, this implies that $\beta > 0$ and the geodesic to consider is $(x^{+}(\cdot, \beta), y(\cdot, \beta))$. The first equation in \eqref{geodesicsagrushin0} implies that for a solution to exist, we must have $\beta \leqslant \kappa/x_1^\alpha$. When that is the case, there are two solutions:
	\begin{equation}
		\label{timeequation}
		\left\{
		\begin{array}{rcl}
			t_1(\beta) & = & \dfrac{1}{\beta} \Biggl(\dfrac{\beta}{\kappa}\Biggr)^{\frac{\alpha - 1}{\alpha}} \arcsin_\alpha \Biggl( \dfrac{x_1 \beta^{1/\alpha}}{\kappa^{1/\alpha}} \Biggr) \\
			t_2(\beta) & = & \dfrac{1}{\beta} \Biggl(\dfrac{\beta}{\kappa}\Biggr)^{\frac{\alpha - 1}{\alpha}} \Biggl[\pi_\alpha -  \arcsin_\alpha \Biggl( \dfrac{x_1 \beta^{1/\alpha}}{\kappa^{1/\alpha}} \Biggr)\Biggr].
		\end{array}
		\right.
	\end{equation}
	The function $t_1(\beta)$ is increasing from $x_1/\kappa$ as $\beta$ goes to 0, to $\pi_\alpha/|\omega|$ when $\beta = \kappa/x_1^{\alpha}$. The function $t_2(\beta)$ is decreasing from $+\infty$ when $\beta$ tends to $0$, to $\pi_\alpha/|\omega|$ when $\beta = \kappa/x_1^\alpha$. We substitute these two into the second equation in \eqref{geodesicsagrushin0} and use the identity $\cos_\alpha^2(x) = 1 - \sin_\alpha^{2 \alpha}(x)$. The assumption $y_1 \geqslant y_0$ enables us to choose the positive sign when taking the square root:
	\begin{equation}
		\label{yequation}
		\left\{
		\begin{array}{rcl}
			y_1(\beta) & = & y_0 + \dfrac{(\kappa/\beta)^{\frac{\alpha + 1}{\alpha}}}{(\alpha + 1)} \Biggl[ \arcsin_\alpha \left( \dfrac{x_1 \beta^{1/\alpha}}{\kappa^{1/\alpha}} \right) - \sqrt{ 1 - \dfrac{x_1^{2 \alpha} \beta^2}{\kappa^2}} \dfrac{x_1 \beta^{1/\alpha}}{\kappa^{1/\alpha}} \Biggr] \\
			y_2(\beta) & = & y_0 + \dfrac{(\kappa/\beta)^{\frac{\alpha + 1}{\alpha}}}{(\alpha + 1)} \Biggl[\pi_\alpha -  \arcsin_\alpha \left( \dfrac{x_1 \beta^{1/\alpha}}{\kappa^{1/\alpha}} \right) + \sqrt{ 1 - \dfrac{x_1^{2 \alpha} \beta^2}{\kappa^2}} \dfrac{x_1 \beta^{1/\alpha}}{\kappa^{1/\alpha}} \Biggr].
		\end{array}
		\right.
	\end{equation}
	
	The function $y_1(\beta)$ is increasing (resp. $y_2(\beta)$ is decreasing) and behaves in the following way. When $\beta$ tends to 0, $y_1$ goes to $y_0$ (resp. $y_2$ goes to $+\infty$) and when $\beta = \kappa/x_1^{\alpha}$, the function $y_1$ (resp. $y_2$) takes the value $y_0 + x_1^{\alpha + 1} \pi_\alpha/[2(\alpha + 1)]$. Therefore, given $x_1 > 0$ and $y_1 \geqslant y_0$, if $y_1 \leqslant y_0 + x_1^{\alpha + 1} \pi_\alpha/[2(\alpha + 1)]$ (resp. $y_1 \geqslant y_0 + x_1^{\alpha + 1} \pi_\alpha/[2(\alpha + 1)]$), we use \eqref{yequation} to deduce the existence of a unique $\beta > 0$ such that $y_1(\beta) = y_1$ (resp. $y_2(\beta) = y_1$) and \eqref{timeequation} provides the unique $t = t_1(\beta) \in \interval[open right]{0}{\pi_\alpha / |\omega|}$ (resp. $t = t_2(\beta)$) such that $(x^{+}(t, \beta), y(t, \beta)) = (x_1, y_1)$.
	
	The geodesic $\gamma$ is consequently minimising before $t = \pi_\alpha / |\omega|$.
\end{proof}

It remains to study the case of a geodesic $\gamma$ starting at a Riemannian point $(x_0, y_0)$, i.e. with $x_0 \neq 0$. We will use an extended Hadamard technique, as described in \cite[Section 13.4]{comprehensivesubr}:

\begin{theorem}[Extended Hadamard technique]
	\label{thm:carthad}
	Let $M$ be an ideal sub-Riemannian manifold and $q_0 \in M$ be a Riemannian point (resp. a singular point). Let $\mathrm{Cut}^*(q_0) \subseteq M$ be the conjectured cut locus and $t^*_{q_0}[\lambda_0] \in \interval{0}{+\infty}$ be the conjectured cut time at $q_0$ for an initial covector $\lambda_0 \in \mathrm{T}^*_{q_0}(M) \cap H^{-1}(1/2)$. 
	
	Set $N$ as the set of covectors in $\mathrm{T}^*_{q_0}(M)$ for which the corresponding geodesics are conjectured to be optimal up to time 1. 
	
	In other words,
	\[
	N := \{t \theta \mid \lambda_0 \in \mathrm{T}^*_{q_0}(M) \cap H^{-1}(1/2) \text{ and } t \in \interval[open right]{0}{t^*_{q_0}[\lambda_0]} \ (\text{resp. }t \in \interval[open left, open right]{0}{t^*_{q_0}[\lambda_0]})\}.
	\]
	Assume that the set $N$ is shown to satisfy the following conditions:
	\begin{enumerate}[label=\normalfont(\roman*)]
		\item $\mathrm{exp}_{q_0}(N) = M \setminus \mathrm{Cut}^*(q_0)$;
		\item The restriction of the sub-Riemannian exponential $\mathrm{exp}_{q_0}|_N$ is a proper map, invertible at every point of $N$;
		\item The set $\mathrm{exp}_{q_0}(N)$ is simply-connected (resp. $\mathrm{exp}_{q_0}|_N$ is a diffeomorphism).
	\end{enumerate}
	Then, $\mathrm{exp}_{q_0}|_N$ is a diffeomorphism and the conjectured cut locus and cut times are the right ones: $\mathrm{Cut}(q_0) = \mathrm{Cut}^*(q_0)$ and $t_{q_0} = t^*_{q_0}$. 
\end{theorem}

\begin{remark}
	The restriction of $\mathrm{T}^*_{q_0}(M)$ to $H^{-1}(1/2)$ results from considering geodesics parametrised by arclength.
\end{remark}

We firstly observe that 
\[
\gamma\left( \dfrac{\pi_\alpha}{|\omega|} \middle| A, \omega, \phi \right) = \gamma\left( \dfrac{\pi_\alpha}{|\omega|} \middle| A, \omega, \pi_\alpha - \phi \right).
\]
This means that the points
\[
\left(- x_0, y_0 + \mathrm{sgn}(\omega) \left(\dfrac{x_0}{\sin_\alpha(\phi)}\right)^{\alpha + 1} \dfrac{\pi_\alpha}{(\alpha + 1)}\right)
\]
are joined from $(x_0, y_0)$ by two distinct geodesics unless $\phi = \pi_\alpha/2$ or $3\pi_\alpha/2$ in which case there is only one.

This leads us to conjecture that the cut time should be $t^*_\mathrm{cut}(u_0, v_0) = \pi_\alpha / |\omega|$ and that the cut locus should be
\[
\mathrm{Cut}^*(x_0, y_0) = \left\{(- x_0, y) \in \mathds{G}_\alpha \mid |y - y_0| \geqslant |x_0|^{\alpha + 1} \dfrac{\pi_\alpha}{(\alpha + 1)}\right\}.
\]

Here, the set of covectors in $\mathrm{T}^*_{q_0}(\mathds{G}_\alpha)$ for which the corresponding geodesics are conjectured to be optimal up to time 1 is 
\begin{equation}
	\label{cotinjec}
	\begin{array}{rcl} 
		N & := & \left\{ t \lambda_0 \mid \lambda_0 \in \mathrm{T}^*_{(x_0, y_0)}(\mathds{G}_\alpha) \cap H^{-1}(1/2), \ t \in \rinterval{0}{t^*_\mathrm{cut}[\lambda_0]} \right\}  \\
		& = &  \left\{ u_0 \mathrm{d}x|_{(x_0, y_0)} + v_0 \mathrm{d}y|_{(x_0, y_0)} \in \mathrm{T}_{(x_0, y_0)}^*(\mathds{G}_\alpha) \mid |\omega| < \pi_\alpha \right\},
	\end{array}
\end{equation}
and thus $\exp_{(x_0, y_0)}(N) = \left\{ (x, y) \in \mathds{G}_\alpha \mid (x, y) \notin \mathrm{Cut}^*(q_0) \right\}$.

Let us show that the equality in \cref{cotinjec} indeed holds. When considering a covector $\lambda_0$ (resp. $\overline{\lambda}_0$), we write $A$ and $\omega$ (resp. $\overline{A}$ and $\overline{\omega}$) for the corresponding coordinates given by \eqref{corresp}. If $v_0$ tends to $0$, then $\omega$ tends to $0$ and $t^*_\mathrm{cut}[\lambda_0] = +\infty$ which implies that covectors with $v_0 = 0$ belong to both sets in \eqref{cotinjec}. We can now assume that $v_0 \neq 0$ (resp. $\overline{v}_0 \neq 0$). If $\overline{\lambda}_0 = t \lambda_0$ is a vector in $N$ for some $t \in \rinterval{0}{t^*_\mathrm{cut}[\lambda_0]}$, then, with the help of \eqref{corresp}, we find that $|\overline{\omega}| = t |\omega|$ and therefore $|\overline{\omega}| < \pi_\alpha$. On the other hand, if $\overline{\lambda_0}$ is a covector such that $|\overline{\omega}| < \pi_\alpha$, we can express it as $\overline{\lambda}_0 = t \lambda_0$ with $t := \overline{A} \overline{\omega} > 0$ and $\lambda_0 := \overline{\lambda}_0/t$. Using \eqref{corresp} again, we deduce that $A \omega = 1$ and thus $\lambda_0 \in H^{-1}(1/2)$. Furthermore, the coefficient $t$ satisfies
\[
0 \leqslant t = |\overline{A}| |\overline{\omega}| = |A| |\overline{\omega}| = \dfrac{|\overline{\omega}|}{|\omega|} < \dfrac{\pi_\alpha}{|\omega|},
\]
since $|\overline{\omega}| < \pi_\alpha$ by hypothesis.

\begin{figure}
	\centering
	\captionsetup[subfigure]{justification=centering}
	\subcaptionbox{$x_0 = 0, \alpha = 1$\label{fig3:a}}{\scalebox{0.465}{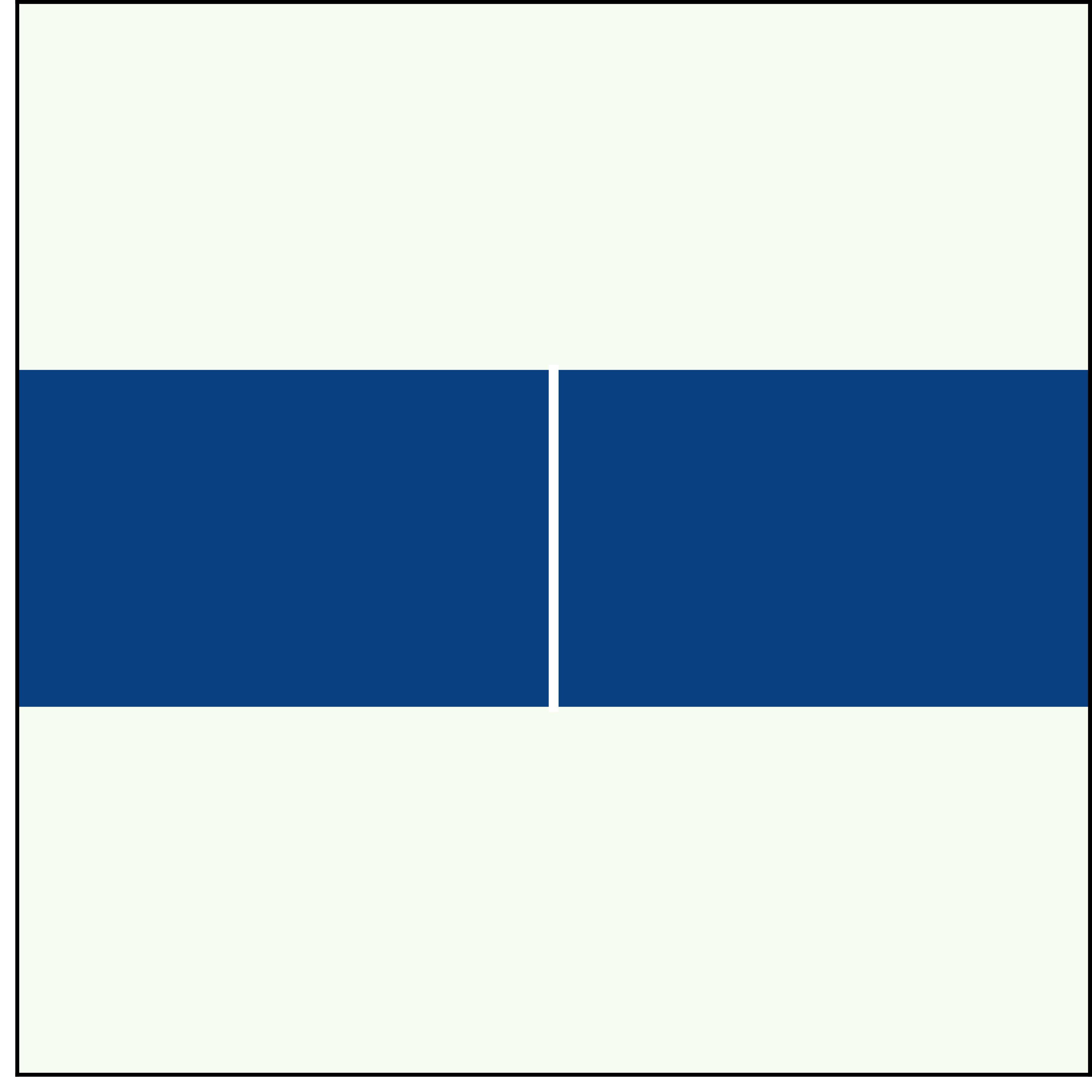}}\hspace{1em}%
	\subcaptionbox{$x_0 = 0, \alpha > 1 $\label{fig3:b}}{\scalebox{0.465}{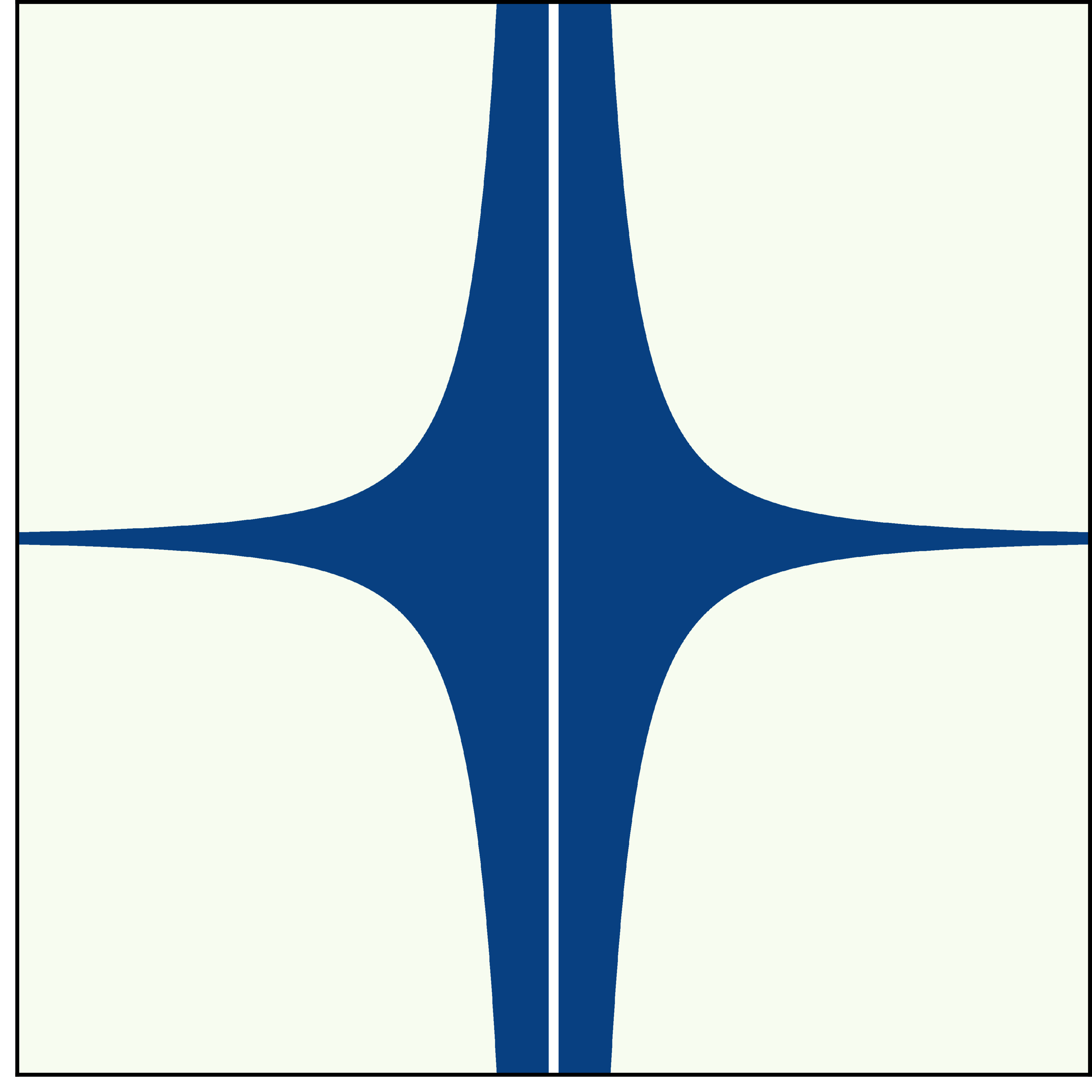}}\hspace{1em}%
	\subcaptionbox{$x_0 \neq 0, \alpha = 1 $\label{fig3:c}}{\scalebox{0.465}{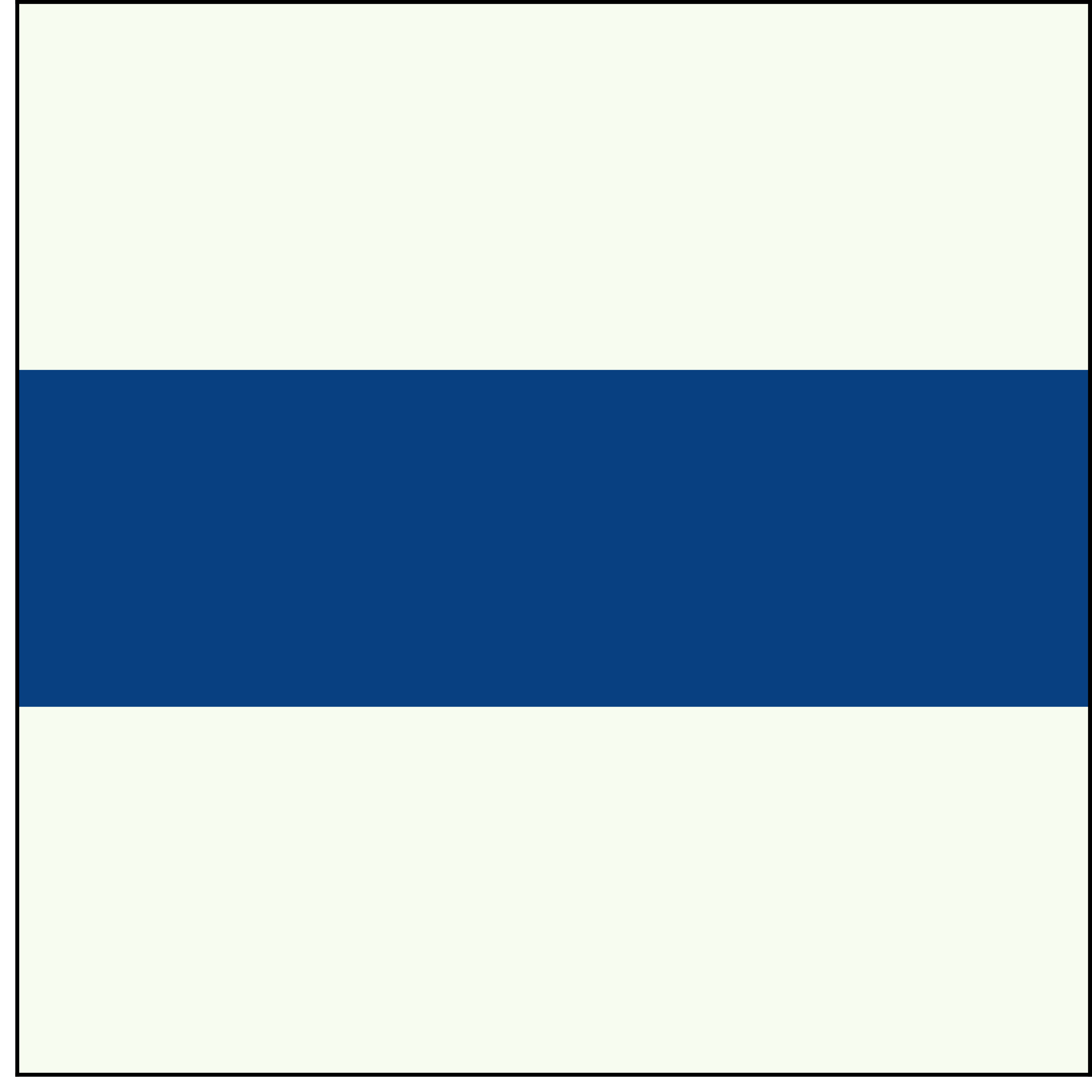}}\hspace{1em}%
	\subcaptionbox{$x_0 \neq 0, \alpha > 1 $\label{fig3:d}}{\scalebox{0.465}{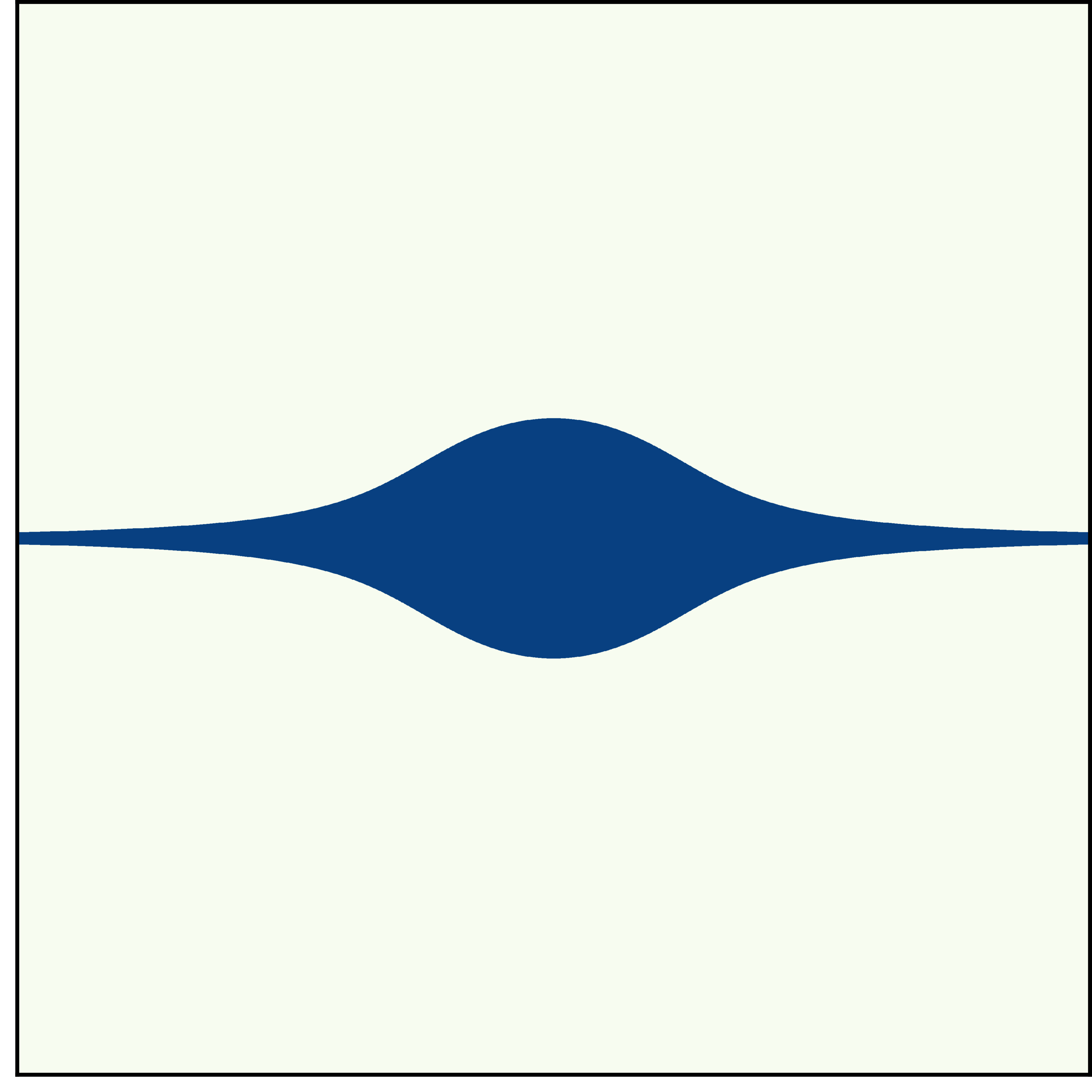}}\hspace{1em}%
	\caption{The cotangent injectivity domain is an open star-shaped region if $x_0 \neq 0$. If $x_0 = 0$, it looks like a star-shaped region but with the starting point and the annihilator of the distribution removed.}
	\label{injradius}
\end{figure}

\begin{remark}
	The set \eqref{cotinjec} corresponds to what is called the \textit{(cotangent) injectivity domain}. If $x_0 = 0$, the cotangent injectivity domain will be as in \eqref{cotinjec} but with $H^{-1}(0)$ being removed, since this time $t \in \interval[open left, open right]{0}{t^*_{q_0}[\lambda_0]}$ by \cref{thm:carthad}. When $\alpha = 1$, the condition defining $N$ reduces to $|v_0| \leqslant \pi$. Geometrically, this is a horizontal strip in the cotangent space. The shape of the cotangent injectivity domain for $\alpha > 1$ is different than when $\alpha = 1$: see \cref{injradius}.
\end{remark}

We know that $A^2 \omega^2 = u_0^2 + v_0^2 x_0^{2 \alpha} = \kappa^2 = 2 H(u_0, v_0)$, where the positive parameter $\kappa > 0$ is the constant speed of the geodesic $\gamma$. We can then parametrise $u_0, v_0$ and the corresponding parameters $A$ and $\omega$ with respect to $t \in \interval{0}{t^*_\mathrm{cut}(u_0, v_0)}$ and $\phi \in \ointerval{0}{2 \pi_\alpha} \setminus \{ \pi_\alpha\}$:
\[
u_0 = \kappa \cos_\alpha(\phi), \ \ v_0 = \kappa \dfrac{\sin_\alpha(\phi)}{x_0} \left|\dfrac{\sin_\alpha(\phi)}{x_0}\right|^{\alpha - 1}, A = \dfrac{x_0}{\sin_\alpha(\phi)} \text{ and } \omega = \kappa \dfrac{\sin_\alpha(\phi)}{x_0}.
\]
The expression of the geodesics from \cref{thm:geodesics} can thus be written as
\begin{equation}
	\label{geodesicsagrushin}
	\left\{
	\begin{array}{rcl}
		x(t, \phi) & = & \dfrac{x_0}{\sin_\alpha(\phi)} \sin_\alpha\left( \kappa \dfrac{\sin_\alpha(\phi)}{x_0} t + \phi \right) \\
		y(t, \phi) & = & y_0 + \dfrac{1}{(\alpha + 1)}\left|\dfrac{x_0}{\sin_\alpha(\phi)}\right|^{\alpha + 1}\left[ \kappa \dfrac{\sin_\alpha(\phi)}{x_0} t + \cos_\alpha(\phi) \sin_\alpha(\phi) \right.\\
		& & \left. - \cos_\alpha\left( \kappa \dfrac{\sin_\alpha(\phi)}{x_0} t + \phi \right) \sin_\alpha\left( \kappa \dfrac{\sin_\alpha(\phi)}{x_0} t + \phi \right) \right].
	\end{array}
	\right.
\end{equation}
In fact, $\phi = 0$ or $\pi_\alpha$ correspond to the geodesic starting at $(x_0, y_0)$ with initial covector $(\kappa, 0)$ and $(-\kappa, 0)$ respectively. In that case, the geodesics are parametrised by
\begin{equation}
	\label{geodesicsagrushin2}
	\left\{
	\begin{array}{rcl}
		x(t, 0) & = & \kappa t + x_0 \\
		x(t, \pi_\alpha) & = & - \kappa t + x_0 \\
		y(t, 0) & = & y_0 \\
		y(t, \pi_\alpha) & = & y_0
	\end{array}
	\right.
\end{equation}
Given a constant speed $\kappa > 0$ and an initial point $p := (x_0, y_0)$ with $x_0 \neq 0$, we can compute the determinant of the differential of the corresponding \textit{exponential map} $(t, \phi) \mapsto (x(t, \phi), y(t, \phi))$:
\begin{equation}
	\label{detx0not0}
	\begin{array}{rcl} 
		D(t, \phi) & = & \dfrac{\kappa}{x_0 \sin_\alpha(\phi)} \left|\dfrac{x_0}{\sin_\alpha(\phi)}\right|^{\alpha + 1} \left[x_0 \sin_\alpha\left( \kappa \dfrac{\sin_\alpha(\phi)}{x_0} t + \phi \right) \cos_\alpha(\phi)\right. \\
		& & \left. - \sin_\alpha(\phi) (x_0 + \kappa t \cos_\alpha(\phi)) \cos_\alpha\left( \kappa \dfrac{\sin_\alpha(\phi)}{x_0} t + \phi \right) \right]. 
	\end{array}
\end{equation}
One can check that $\lim\limits_{\phi \to 0} D(t, \phi) = \lim\limits_{\phi \to \pi_\alpha} D(t, \phi) = 0$ unless $\alpha = 1$, in which case we have 
\[
\lim\limits_{\phi \to 0} D(t, \phi) = \dfrac{\kappa^4 t}{3 x_0^5}\left(\kappa^2 t^2 + 3 \kappa t x_0 + 3 x_0^2 \right) 
\]
and
\[
\lim\limits_{\phi \to 0} D(t, \phi) = \dfrac{\kappa^4 t}{3 x_0^5}\left(\kappa^2 t^2 - 3 \kappa t x_0 + 3 x_0^2 \right).
\]

We now claim that the exponential map has no singularities before $t = \pi_\alpha / |\omega|$. Indeed, we observe firstly that $D(0, \phi)$ vanishes for every $\phi$. Secondly, with the help of the derivative of $D$ with respect to $t$;
\begin{equation*}
	\begin{array}{rcl} 
		\partial_t D(t, \phi) & = & \alpha \dfrac{\kappa^2}{x_0^2} (x_0 + \kappa t \cos_\alpha(\phi)) \left|\dfrac{x_0}{\sin_\alpha(\phi)}\right|^{\alpha + 1} \sin_\alpha(\phi) \\
		& & \times \sin^{2 (\alpha - 1)}_\alpha\left( \kappa \dfrac{\sin_\alpha(\phi)}{x_0} t + \phi \right) \sin_\alpha\left( \kappa \dfrac{\sin_\alpha(\phi)}{x_0} t + \phi \right), 
	\end{array}
\end{equation*}
we see that $\partial_t D(t, \phi) = 0$ if and only if 
\[
t = - \dfrac{x_0}{\kappa \cos_\alpha(\phi)} \text{ or } t = \dfrac{x_0}{\kappa \sin_\alpha(\phi)} ( l \pi_\alpha - \phi), \ l \in \mathds{Z}.
\]
The former is a local minimum that is positive while the later is a local maximum that is also positive. Thirdly, we observe that
\[
D(t^*_\mathrm{cut}, \phi) = \kappa \dfrac{\pi_\alpha}{\sin_\alpha(\phi)} \left|\dfrac{x_0}{\sin_\alpha(\phi)}\right|^{\alpha + 1} \cos^2_\alpha(\phi),
\]
which is zero if and only if $\phi = \pi_\alpha / 2$ or $3 \pi_\alpha/2$. So, the function $D$ is never zero on $\ointerval{0}{t^*_\mathrm{cut}}$ and the exponential map is invertible at every point of $N$.

Finally, we need to make some topological considerations in order to conclude. Consider the set $N$ for which the corresponding geodesics are conjectured to be optimal up to time 1 and its image under the sub-Riemannian exponential map at $(x_0, y_0)$. The map $\exp_{(x_0, y_0)} : N \to \exp(N)$ is proper: if a sequence of points $(u_i, v_i) \in N$ escapes to infinity, we must have $u_i \to \pm\infty$ and therefore $\exp_{(x_0, y_0)}(u_i, v_i)$ will also escape to infinity. Therefore, $\exp|_N$ is indeed proper, its differential is not singular at any point and furthermore $\exp(N)$ is simply connected. 

We can conclude that $\exp$ is a diffeomorphism and the extended Hadamard technique (\cref{thm:carthad}) implies that the conjectured cut loci and time are thus the true ones.

\begin{figure}
	\centering
	\def\svgwidth{\columnwidth}
	\scalebox{0.494}{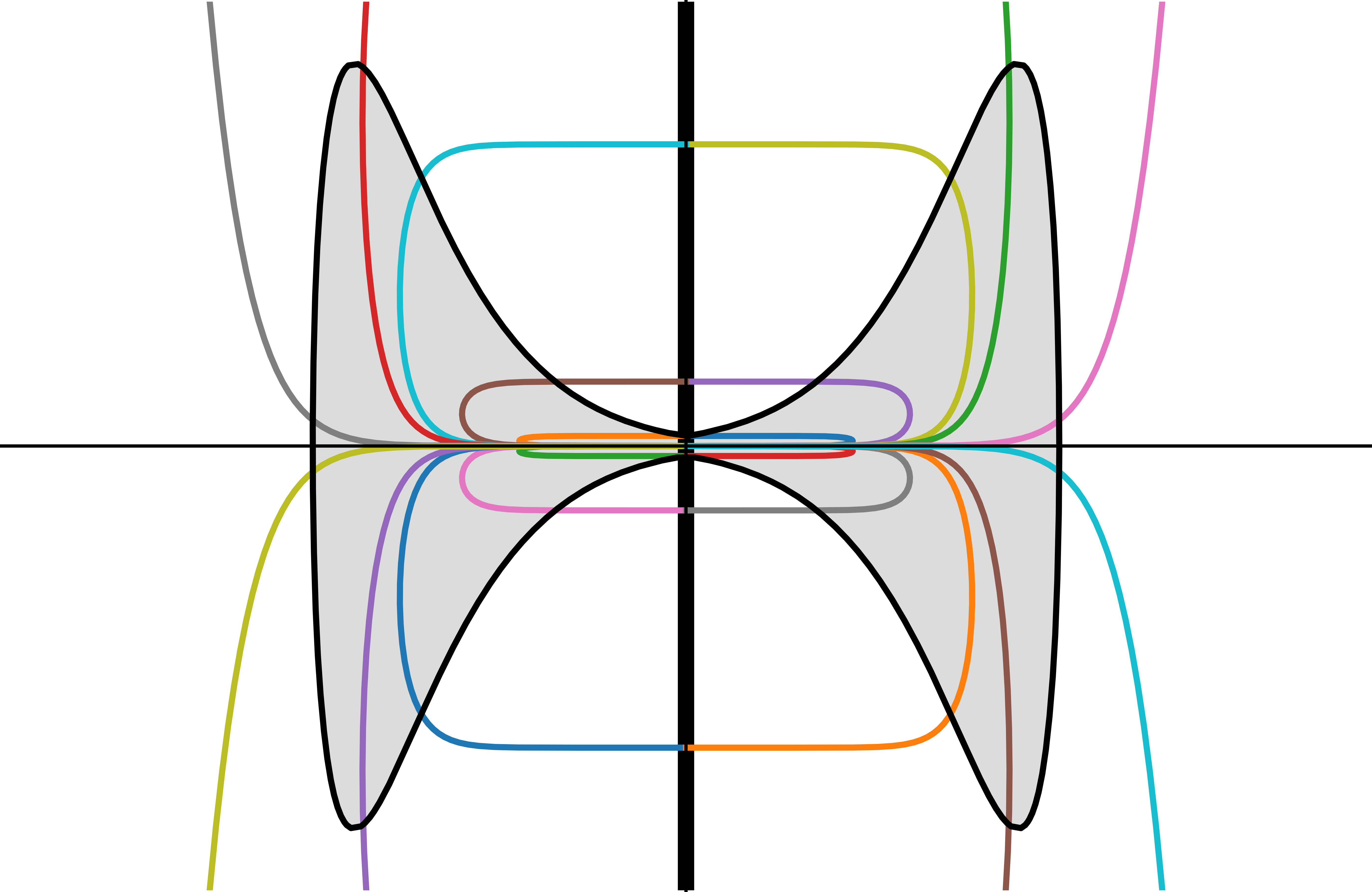}
	\scalebox{0.665}{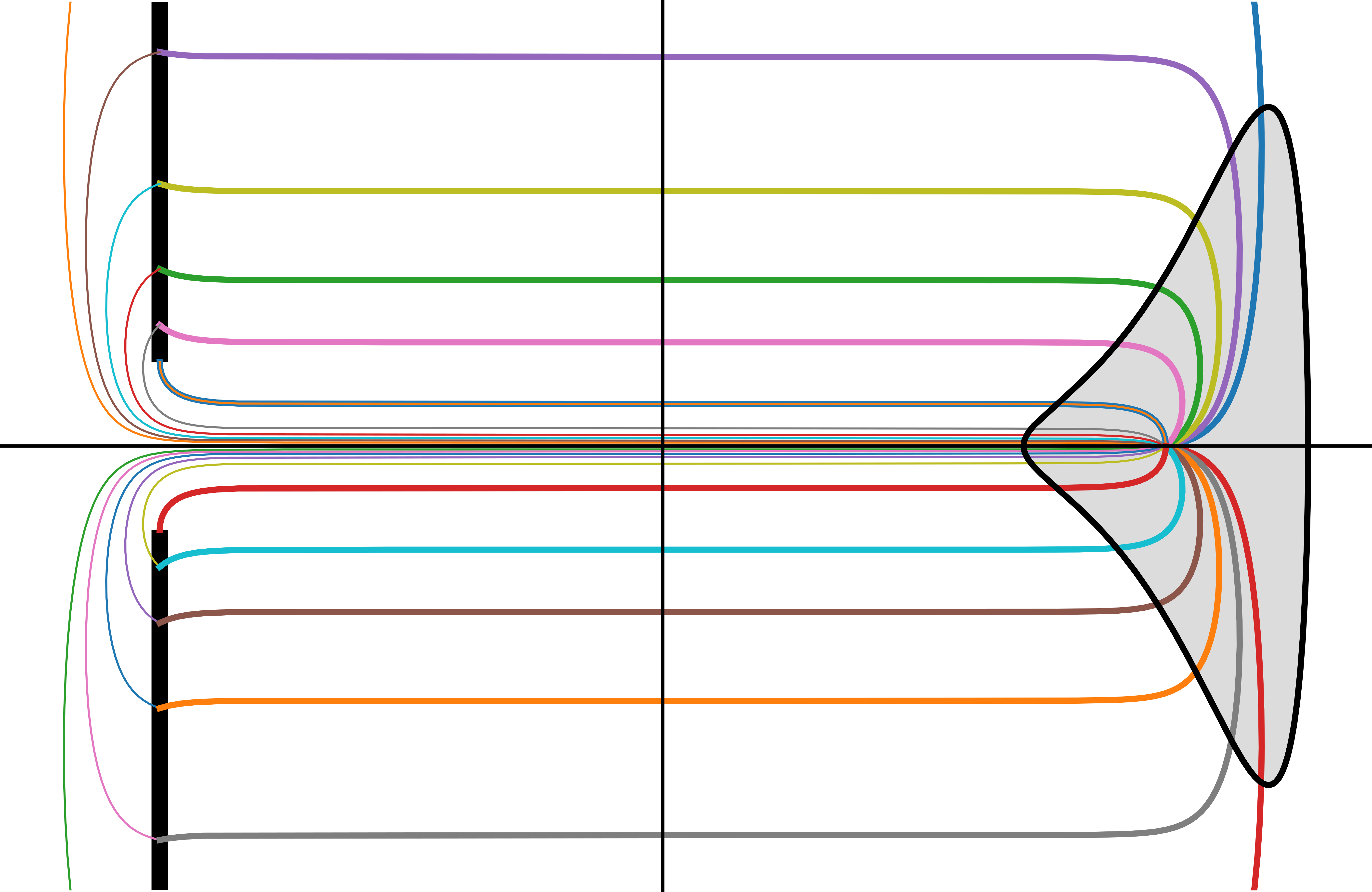}
	\caption{Geometry of $\mathds{G}_\alpha$}
	Illustration of the geodesics of the $\alpha$-Grushin plane from a singular point (on the left) and from a Riemannian point (on the right). The shaded area represents a ball around the starting point and the thick line is the cut locus.
	\label{fig:geodesics}
\end{figure}

To summarise the findings of this section, we have proved the following result:

\begin{theorem}
	\label{cutlocus}
	Let $\alpha \geqslant 1$ and $\gamma(t) = (x(t), y(t))$ be a geodesic of $\mathds{G}_\alpha$ with initial value $\gamma(0) = (x_0, y_0)$ and initial covector $u_0 \mathrm{d}x|_{(x_0, y_0)} + v_0 \mathrm{d}y|_{(x_0, y_0)}$, as described in \cref{thm:geodesics}. 
	
	If $v_0 = 0$, there are no singularities along $\gamma$,
	\[
	t_\mathrm{cut}[\gamma] = +\infty \text{ and } \mathrm{Cut}(x_0,y_0) = \emptyset. 
	\] 
	
	If $v_0 \neq 0$, then the cut time is
	\[t_\mathrm{cut}[\gamma] = \dfrac{\pi_\alpha}{|\omega|},\]
	while the cut locus is 
	\[\mathrm{Cut}(x_0,y_0) = \left\{(- x_0, y) \in \mathds{G}_\alpha \mid |y - y_0| \geqslant |x_0|^{\alpha + 1} \dfrac{\pi_\alpha}{(\alpha + 1)}\right\}.
	\]
\end{theorem}
The cut loci and geodesics of $\mathds{G}_\alpha$ are illustrated in  \cref{fig:geodesics}. With this in mind, we now turn to the analysis of the distortion coefficients of the $\alpha$-Grushin plane.

\section{Distortion coefficients of the \texorpdfstring{$\alpha$}{alpha}-Grushin plane} \label{distortionagrushin}

\subsection{Computation of the distortion coefficients} \label{distcoeff}

We present here our main result: an explicit computation of the distortion coefficient of $\mathds{G}_\alpha$. To this aim, we use the techniques established by Balogh, Kristály and Sipos in \cite{baloghheis} and generalised by Barilari and Rizzi in \cite{subriemint}. In the latter, the authors prove interpolation inequalities of optimal transport for ideal sub-Riemannian manifolds. They are expressed in terms of the distortion coefficients for which the expression is obtained through a fine analysis of sub-Riemannian Jacobi fields.

\begin{theorem} \label{thm:distcoef2}
	Let $q, q_0 \in \mathds{G}_\alpha$ such that $q \notin \mathrm{Cut}(q_0)$. Assume that $q$ and $q_0$ do not lie on the same horizontal line. Under the correspondence of \cref{thm:geodesics} and the relations \eqref{corresp}, we have
	\[
	\beta_t(q, q_0) = \dfrac{\mathrm{J}(t, A, \omega, \phi)}{\mathrm{J}(1, A, \omega, \phi)} \text{  for all } t \in \interval{0}{1},
	\]
	where
	\begin{equation}
		\label{distAphiom}
		\mathrm{J}(t, A, \omega, \phi) = t \Bigl[\sin_\alpha(\omega t + \phi) \cos_\alpha(\phi) - \cos_\alpha(\omega t + \phi) \bigl(\sin_\alpha(\phi) + \omega t \cos_\alpha(\phi)\bigr)\Bigr].
	\end{equation}
\end{theorem}

\begin{remark}
	We consider geodesics parametrised by constant speed on $\interval{0}{1}$. Consequently, since \cref{cutlocus} states that $t_\mathrm{cut} = \pi_\alpha/|\omega|$, we always have $|\omega| \leqslant \pi_\alpha$ when $q \notin \mathrm{Cut}(q_0)$.
\end{remark}

\begin{proof}
	We let $\lambda_0 = u_0 \mathrm{d}x|_{q_0} + v_0 \mathrm{d}y|_{q_0} \in \mathrm{T}^*_{q_0}(\mathds{G}_\alpha)$ be the covector corresponding to the unique minimising geodesic joining $q_0 = (x_0, y_0)$ to $q = (x, y)$ in $\mathds{G}_\alpha$. The assumption that $q$ and $q_0$ do not lie on the same horizontal line means that $v_0 \neq 0$.
	
	By choosing the global Darboux frame induced by the sections of $\mathrm{T}(\mathrm{T}^*(\mathds{G}_\alpha))$; $E_1 = \partial_u, E_2 = \partial_v, F_1 = \partial_x, F_2 = \partial_y$, Lemma 44 in \cite{subriemint} yields that $\beta_t(q,q_0) = J(t)/J(1)$	where the function $J$ is the determinant of the exponential map $(u, v) \to \mathrm{exp}_{(x_0, y_0)}(u, v)$ in these coordinates, computed at $(u_0, v_0)$. 
	
	Taking the derivatives of \eqref{geodesics}, we find
	\begin{equation*}
		\begin{array}{rcl}
			x_{u_0} (t) & = & A_{u_0} \sin_\alpha(\omega t + \phi) + A \left(\omega_{u_0} t + \phi_{u_0} \right) \cos_\alpha(\omega t + \phi); \\
			x_{v_0} (t) & = & A_{v_0} \sin_\alpha(\omega t + \phi) + A \left(\omega_{v_0} t + \phi_{v_0} \right) \cos_\alpha(\omega t + \phi),
		\end{array}
	\end{equation*}
and
\begin{equation*}
	\begin{array}{rcl}
		y_{u_0}(t) & = & \dfrac{v_0 A^{2 (\alpha - 1)} A}{(\alpha + 1) \omega^2} \Bigl[ \alpha \omega t \bigl(2 A_{u_0} \omega + A  \omega_{u_0} \bigr) \\
		& & \qquad + A \omega (\alpha + 1) (\phi_{u_0} \cos_\alpha^2(\phi)- \left(\omega_{u_0} t + \phi_{u_0} \right) \cos_\alpha^2(\omega t + \phi)) \\
		& & \qquad + \bigl(2 \alpha A_{u_0} \omega - A \omega_{u_0} \bigr) \\
		& & \qquad \qquad \times \bigl(\sin_\alpha(\phi) \cos_\alpha(\phi) - \sin_\alpha(\omega t + \phi) \cos_\alpha(\omega t + \phi) \bigr) \Bigr]; \\
		y_{v_0}(t) & = & \dfrac{v_0 A^{2 (\alpha - 1)} A}{(\alpha + 1) \omega^2} \Bigl[ \omega t \bigl(2 \alpha v_0 \omega A_{v_0} + A (\omega + \alpha v_0 \omega_{v_0})\bigr) \\
		& & \qquad + A \omega v_0 (\alpha + 1) (\phi_{u_0} \cos_\alpha^2(\phi)- \bigl(\omega_{v_0} t + \phi_{v_0} \bigr) \cos_\alpha^2(\omega t + \phi)) \\
		& & \qquad + \bigl(2 \alpha \omega v_0 A_{v_0} + A (\omega - v_0 \omega_{v_0}) \bigr) \\
		& & \qquad \qquad \times \bigl(\sin_\alpha(\phi) \cos_\alpha(\phi) - \sin_\alpha(\omega t + \phi) \cos_\alpha(\omega t + \phi) \bigr) \Bigr].
	\end{array}
\end{equation*}
	To make things clearer, set $[f, g] := f_{u_0} g_{v_0} - f_{v_0} g_{u_0}$ and we obtain
	\begin{equation*}
	\begin{array}{rl}
		[x, y](t) = & \hspace{-7pt} \dfrac{A^{2 \alpha}}{(\alpha + 1) \omega^2} \Bigl[\sin_\alpha^2(\omega t + \phi) \cos_\alpha(\phi) \bigl( v [A, \omega] - \omega A_{u_0} \bigr) \\
		& + \alpha v_0 \omega \sin_\alpha^{2 \alpha}(\omega t + \phi) \sin_\alpha(\omega t + \phi) \bigl([A, \omega] t + [A, \phi] \bigr) \\
		& + \sin_\alpha(\omega t + \phi) \Bigl( \sin_\alpha(\phi) \cos_\alpha(\phi) (\omega A_{u_0} - v_0 [A, \omega]) \\
		& \quad - \alpha v_0 \omega \sin_\alpha^{2 \alpha}(\phi) [A, \phi] + \omega (\omega A_{u_0} t + v_0 \cos_\alpha^2(\phi) [A, \phi] ) \Bigr) \\
		& + \sin_\alpha(\omega t + \phi) \cos^2_\alpha(\omega t + \phi) \Bigl( t \bigl( (2 \alpha - 1) v_0 \omega [A, \omega] - A \omega \omega_{u_0}\bigr) \\
		& \quad + (2 \alpha - 1) v_0 \omega [A, \phi]  + A v_0[\phi, \omega] - A \omega \phi_{u_0} \Bigr) \\
		& + \cos_\alpha(\omega t + \phi) \Bigl( \sin_\alpha(\phi) \cos_\alpha(\phi) (A \omega \phi_{u_0} + 2 \alpha v_0 \omega [\phi, A] + A v_0 [\omega, \phi]) \\
		& \quad + \omega^2 t^2 (2 \alpha v_0 [\omega, A] + A \omega_{u_0})  \\
		& \quad - \omega t ( \sin_\alpha(\phi) \cos_\alpha(\phi) (2 \alpha v_0 [A, \omega]) - A \omega_{u_0})  \\
		& \quad + \omega (2 \alpha v_0 [A, \phi]) - A \phi_{u_0} + \alpha v_0 A \sin_\alpha^{2 \alpha}(\phi) [\omega, \phi] + \alpha v_0 \cos_\alpha^2(\phi) [\phi, \omega] \Bigr)\Bigr].
	\end{array}
	\end{equation*} 
	Using the identities $\sin_\alpha^{2 \alpha}(x) + \cos_\alpha^2(x) = 1$ and \eqref{der}, we find that $[\omega, \phi] =  \frac{\sin_\alpha(\phi)}{\alpha v_0 A}$, $[A, \omega] = \frac{\cos_\alpha(\phi)}{\alpha v_0}$, and $[A, \phi] = 0$.
	Consequently, we have
	\begin{equation}
		\label{jacobianexp}
		\begin{array}{rl}
			[x, y](t) = & \hspace{-5pt} t \dfrac{A^{2 \alpha}}{\alpha \omega} \Bigl[\sin_\alpha(\omega t + \phi) \cos_\alpha(\phi) - \cos_\alpha(\omega t + \phi) \bigl(\sin_\alpha(\phi) + \omega t \cos_\alpha(\phi)\bigr)\Bigr].
		\end{array}
	\end{equation}
	We finally obtain the desired expression by performing $\beta_t(q, q_0) = [x, y](t)\big/[x, y](1)$.
\end{proof}

We can transform \eqref{distAphiom} from the set of coordinates $A, \omega$ and $\phi$ to $x_0, u_0$ and $v_0$ via the identities \eqref{corresp}. This leads to \eqref{eq:jacdetxuv} and concludes the proof of \cref{thm:distcoef}.

It is interesting to look at the limit of \eqref{eq:jacdetxuv} when $\alpha$ tends to $1$. In this case, the $\alpha$-Grushin plane is the traditional Grushin plane while $\sin_\alpha$ and $\cos_\alpha$ are the usual sine and cosine functions. The formula \eqref{jacobianexp} simplifies to
	\begin{align*}
	[x, y](t) & = t \dfrac{A^{2}}{\omega} \Bigl[ \sin(\omega t  + \phi) \cos(\phi) - \cos(\omega t + \phi) \bigl(\sin(\phi) + \omega t \cos(\phi) \bigr) \Bigr] \\
	& = t \dfrac{(u_0^2 + t u_0 v_0^2 x_0 + v_0^2 x_0^2) \sin(t v_0) - t u_0^2 v_0 \cos(t v_0)}{v_0^3},
	\end{align*}
and thus, we find what was already established in \cite[Proposition 61]{subriemint}: the distortion coefficients of the usual Grushin plane are
\[
\beta_t(q, q_0) = t \dfrac{(u_0^2 + t u_0 v_0^2 x_0 + v_0^2 x_0^2) \sin(t v_0) - t u_0^2 v_0 \cos(t v_0)}{(u_0^2 + u_0 v_0^2 x_0 + v_0^2 x_0^2) \sin(v_0) - u_0^2 v_0 \cos(v_0)} \text{,  for all } t \in \interval{0}{1}.
\]

We now want to investigate the behaviour of $\beta_t(q_0, q)$ when $q_0$ and $q$ do lie on the same horizontal line, that is to say, when $v_0$ tends to $0$. 

\begin{proposition} \label{prop:v0}
	In the same setting of \cref{thm:distcoef2}, when two points $q$ and $q_0$ of $\mathds{G}_\alpha$ are on the same horizontal line, we have
	\[ 
	\beta_t(q_0,q) = t \dfrac{ (u_0 t + x_0)^{2 \alpha}(u_0 t + x_0) - x_0^{2 \alpha}x_0 }{(u_0 + x_0)^{2 \alpha}(u_0 + x_0) - x_0^{2 \alpha}x_0}, \text{ for all } t \in \interval{0}{1}. 
	\]	
\end{proposition}

\begin{remark}
	Considering $q_0$ and $q$ on the same horizontal line corresponds to starting from $q_0$ with an initial covector $u_0 \mathrm{d}x|_{q_0} + v_0 \mathrm{d}y|_{q_0}$ such that $v_0 = 0$. By continuity with respect to initial conditions, the distortion coefficient in \cref{prop:v0} is the limit of \eqref{eq:jacdetxuv} when $v_0$ tends to $0$. In particular, the parameter $u_0$ cannot vanish. Indeed, we would otherwise have a trivial (constant) geodesic since $v_0 = 0$. Furthermore, this implies that the denominator in the expression above is also never vanishing.
\end{remark}

\begin{proof}
	We aim to perform $\lim_{v_0 \to 0} \mathrm{J}(t)\big/\mathrm{J}(1)$, where $\mathrm{J}$ is defined by \eqref{eq:jacdetxuv}. We already know from \cref{thm:geodesics} that $\lim_{v_0 \to 0} u(t) = u_0$ and $\lim_{v_0 \to 0} x(t) = u_0 t + x_0$. Let us make the following preliminary calculations:
	\begin{align*}
		u_{v_0}(t) & = - \alpha A \omega (\omega_{v_0} t + \phi_{v_0}) \sin^{2 (\alpha - 1) }_\alpha(\omega t + \phi) \sin_\alpha(\omega t + \phi) \\
		& \qquad + \cos_\alpha(\omega t + \phi) (\omega A_{v_0} + A \omega_{v_0}) \\
		& = \dfrac{A \omega}{\alpha} \Bigl[ \sin^{2 (\alpha - 1) }_\alpha(\omega t + \phi) \sin_\alpha(\omega t + \phi) \\
		& \qquad \times \Bigl( \omega t ((\alpha - 1) \cos^2_\alpha(\phi) - \alpha) - \sin_\alpha(\phi) \cos_\alpha(\phi) \Bigr) \cos_\alpha(\omega t + \phi) \\
		& \qquad \qquad + ( 1 - \cos^2_\alpha(\phi)) \Bigr] \\
		& = v_0 \dfrac{x_0^{2 \alpha} \cdot u(t) - \left[ t (u_0^2 + \alpha v_0^2 x_0^{2 \alpha}) + u_0 x_0\right] \cdot x(t)^{2 (\alpha - 1)} x(t)}{(u_0^2 + v_0^2 x_0^{2 \alpha})},
	\end{align*}
	and also
	\begin{align*}
		x_{v_0}(t) & = A_{v_0} \sin_\alpha(\omega t + \phi) + A \left(\omega_{v_0} t + \phi_{v_0} \right) \cos_\alpha(\omega t + \phi) \\
		& = \dfrac{A}{\alpha v_0} \Bigl[ \cos_\alpha(\omega t + \phi) \Bigl(\sin_\alpha(\phi) \cos_\alpha(\phi) + \omega t \bigl(\alpha - (\alpha - 1)\cos^2_\alpha(\phi)\bigr) \Bigr) \\
		& \qquad - \sin_\alpha(\omega t + \phi) \cos^2_\alpha(\phi) \Bigr] \\
		& = \dfrac{\left[ t (u_0^2 + \alpha v_0^2 x_0^{2 \alpha}) + u_0 x_0\right] \cdot u(t) + u_0^2 \cdot x(t)}{\alpha v_0 (u_0^2 + v_0^2 x_0^{2 \alpha})}.
	\end{align*}
	Since simply replacing $v_0$ with 0 in $\beta_t(q_0,q)$ will lead to $0/0$, we use L'Hôpital's rule as many times as needed, and we find:
	\begin{equation*}
		\begin{split}
			\beta_t(q_0,q) & = \lim\limits_{v_0 \to 0} \dfrac{\mathrm{J}(t, x_0, u_0, v_0)}{\mathrm{J}(1, x_0, u_0, v_0)} = \lim\limits_{v_0 \to 0} \dfrac{\partial_{v_0} \mathrm{J}(t, x_0, u_0, v_0)}{\partial_{v_0} \mathrm{J}(1, x_0, u_0, v_0)} \\
			& = \lim\limits_{v_0 \to 0} \dfrac{\partial^2_{v_0} \mathrm{J}(t, x_0, u_0, v_0)}{\partial^2_{v_0} \mathrm{J}(1, x_0, u_0, v_0)} = t \dfrac{ (u_0 t + x_0)^{2 \alpha}(u_0 t + x_0) - x_0^{2 \alpha}x_0 }{(u_0 + x_0)^{2 \alpha}(u_0 + x_0) - x_0^{2 \alpha}x_0}.
		\end{split}
	\end{equation*}
\end{proof}

\subsection{Relevant curvature-dimension estimates} \label{mcpagrushin}

Now that we have the expressions for the distortion coefficients, we would like to find appropriate bounds for them. In \cite{juillet2020subriemanniann}, Juillet proved that a sub-Riemannian manifold never satisfies the $\mathrm{CD}(K,N)$ condition when $\mathrm{rank}(\mathcal{D}_p) < \mathrm{dim} (M)$ for all $p \in M$. This result does not apply directly to $\alpha$-Grushin as its distribution has full rank away from the singular set. However, a variant of the technique \cite{juilletdisprove} presented in \cite{juilletthesis} is valid here and we can still conclude that $\mathds{G}_\alpha$ does not satisfy the $\mathrm{CD}$ condition. There is a possibility that the weaker measure contraction property $\mathrm{MCP}(K,N)$ can hold for the $\alpha$-Grushin plane.

In particular, the traditional Grushin plane, equivalent to $\mathds{G}_\alpha$ when $\alpha = 1$, is $\mathrm{MCP}(K, N)$ if and only if $N \geqslant 5$ and $K \leqslant 0$. We expect the $\alpha$-Grushin plane to satisfy the $\mathrm{MCP}$ property for a minimal value of $N$ that would depend on $\alpha$. According to \cref{equivMCP}, the related bound on the distortion coefficients should be of the form $\beta_t(q_0,q) \geqslant t^N$.

In this section, we provide a bound in the case where $q_0$ and $q$ lie on the same horizontal line and when $q_0$ is a Grushin point. In what follows, we will still parametrise the geodesics of the $\alpha$-Grushin plane by constant speed and on the interval $\interval{0}{1}$.

For $\alpha \geqslant 1$, let $m_\alpha \in \interval{-3}{-2}$ be the unique non-zero solution of
\[
(m + 1)^{2 \alpha}(m + 1) - \left( (2 \alpha + 1) m + 1 \right) = 0.
\]
If $\alpha = 1$, the value of the root is $m = -3$.

\begin{proposition}
	\label{prop:distcoefhor}
	Let $q_0 := (x_0, y_0) \in \mathds{G}_\alpha$ with $x_0 \neq 0$ and $q \in \mathds{G}_\alpha$ lying on the same horizontal line. We have that
	\[
	\beta_t(q_0,q) \geqslant t^N \text{ for all } t \in \interval{0}{1}
	\]
	if and only if 
	\[
	N \geqslant 2 \left[ \dfrac{(\alpha + 1) m_\alpha + 1}{m_\alpha + 1} \right].
	\]
\end{proposition}

\begin{proof}
	We are looking for the optimal $N \in \interval{1}{+\infty}$ such that
	\begin{equation}
	\label{eqn:mcpv0}
	\dfrac{ (u_0 t + x_0)^{2 \alpha}(u_0 t + x_0) - x_0^{2 \alpha}x_0 }{(u_0 + x_0)^{2 \alpha}(u_0 + x_0) - x_0^{2 \alpha}x_0} \geqslant t^{N - 1}
	\end{equation}
	for all $ t \in [0,1]$ and $x_0, u_0 \in \mathds{R}$. The function $f_{x_0}(z) := (z + x_0)^{2 \alpha}(z + x_0) - x_0^{2 \alpha}x_0$ is positive (resp. negative) when $z > 0$ (resp. $z < 0$) and $f_{x_0}(0) = 0$. Therefore, the left-hand side of \eqref{eqn:mcpv0} is always non-negative.
	If we take the logarithm of the above, we find that the inequality is equivalent to
	\[
	\int_{u_0 t}^{u_0} \dfrac{\mathrm{d}}{\mathrm{d}z} \log \big| f_{x_0}(z) \big| \mathrm{d} z \leqslant (N - 1) \int_{u_0 t}^{u_0} \dfrac{\mathrm{d}}{\mathrm{d}z} \log |z| \mathrm{d} z.
	\]
	Since this must hold for every $t \in \interval{0}{1}$, it is equivalent to the same inequality for the integrands:
	\[
	\pm \dfrac{(2 \alpha + 1) (z + x_0)^{2 \alpha}}{(z + x_0)^{2 \alpha}(z + x_0) - x_0^{2 \alpha}x_0} \leqslant \pm (N - 1) \dfrac{1}{z}, \text{ when } \pm z > 0.
	\]
	Consequently, \cref{eqn:mcpv0} is equivalent to
	\[
	N \geqslant z\dfrac{(2 \alpha + 1) (z + x_0)^{2 \alpha}}{(z + x_0)^{2 \alpha}(z + x_0) - x_0^{2 \alpha}x_0} + 1.
	\]
	for all $z \in \mathds{R}$ and all $x_0 \in \mathds{R} \setminus \{0\}$. When $z \to 0$, we find that since $x_0 \neq 0$,
	\[
	\dfrac{(2 \alpha + 1) (z + x_0)^{2 \alpha}}{(z + x_0)^{2 \alpha}(z + x_0) - x_0^{2 \alpha}x_0} \rightarrow 1.
	\]
	We are therefore looking for the global maximum of the map
	\begin{align*}
		f \colon \big(\mathds{R} \setminus \{0\}\big) \times \mathds{R} & \longrightarrow \mathds{R} \\
		(x, y) & \longmapsto x\dfrac{(2 \alpha + 1) (x + y)^{2 \alpha}}{(x + y)^{2 \alpha}(x + y) - y^{2 \alpha}y} + 1.
	\end{align*}
	We use polar coordinates: for $r > 0$ and $\theta \in \interval[open right]{0}{2 \pi} \setminus \{\pi/2, -3 \pi/2\}$, we have
	\[
	f(r \cos(\theta), r \sin(\theta)) = \dfrac{(2 \alpha + 1) \cos(\theta) \bigl[\cos(\theta) + \sin(\theta) \bigr]^{2 \alpha}}{\bigl(\cos(\theta) + \sin(\theta)\bigr)\bigl[\cos(\theta) + \sin(\theta) \bigr]^{2 \alpha} - \sin(\theta)^{2 \alpha} \sin(\theta)} + 1,
	\]
	which does not depend on $r$. In particular, the limit of $f$ when $(x, y) \to (0, 0)$ does not exist.
	
	Firstly, let us compute the critical points of $\theta \mapsto f(r \cos(\theta), r \sin(\theta))$. We find that $\frac{\partial}{\partial\theta} f(r \cos(\theta), r \sin(\theta))$ is given by
	\begin{align*}
		& \dfrac{(2 \alpha + 1)\bigl(\cos(\theta) + \sin(\theta)\bigr)\bigl[\cos(\theta) + \sin(\theta)\bigr]^{2(\alpha - 1)}}{\Bigl[\bigl(\cos(\theta) + \sin(\theta)\bigr)\bigl[\cos(\theta) + \sin(\theta) \bigr]^{2 \alpha} - \sin(\theta)^{2 \alpha} \sin(\theta)\Bigr]^2} \\
		& \qquad \times \Bigl[\sin^{2\alpha}(\theta)\bigl((2 \alpha + 1)\cos(\theta) + \sin(\theta)\bigr) \\
		& \qquad \qquad \qquad - \bigl(\cos(\theta) + \sin(\theta)\bigr)\bigl[\cos(\theta) + \sin(\theta) \bigr]^{2 \alpha} \Bigr],
	\end{align*}
	which vanishes when $\cos(\theta) + \sin(\theta) = 0$, i.e. $\theta = 3\pi/4, 7\pi/4$, or when
	\[
	\sin^{2\alpha}(\theta)\left((2 \alpha + 1)\cos(\theta) + \sin(\theta)\right) = \left(\cos(\theta) + \sin(\theta)\right)\left[\cos(\theta) + \sin(\theta) \right]^{2 \alpha}.
	\]
	In the first case, we simply get $f(r \cos(\theta), r \sin(\theta)) = 1$.  The second case implies that $\sin(\theta) \neq 0$, and thus, setting $m = \cot(\theta)$, we obtain
	\begin{equation}
		\label{mequation}
		(m + 1)^{2 \alpha}(m + 1) - \left( (2 \alpha + 1) m + 1 \right) = 0.
	\end{equation}
	\cref{mequation} has two roots: $m = 0$, which we reject since it corresponds to $\theta = \pi/2, 3 \pi/2$ and another root in the interval $\interval{-3}{-2}$, denoted by $m_\alpha$. With a second derivative test, it is easy to see that the $\theta \in \interval[open right]{0}{2 \pi} \setminus \{\pi/2, 3 \pi/2\}$ such that $\theta = \cot^{-1}(m_\alpha)$ gives the local maximums of $\theta \mapsto f(r \cos(\theta), r \sin(\theta))$: at these points, we have
	\begin{align*}
		f(r \cos(\theta), r \sin(\theta)) & = \dfrac{(2 \alpha + 1) \cos(\theta) \left((2 \alpha + 1)\cos(\theta) + \sin(\theta)\right)}{(\cos(\theta) + \sin(\theta))\left[\left((2 \alpha + 1)\cos(\theta) + \sin(\theta)\right) - \sin(\theta)\right]} + 1 \\
		& = \dfrac{(m_\alpha + 1)^{2 \alpha} (2 (\alpha + 1) m_\alpha + 1) - 1}{(m_\alpha + 1)^{2 \alpha}(m_\alpha + 1) - 1} \\
		& = 2 \left[ \dfrac{(\alpha + 1) m_\alpha + 1}{m_\alpha + 1} \right].
	\end{align*}	
	They are in fact global maximums because $f(r \cos(\theta), r \sin(\theta)) \to 1$ when $\theta \to \pi/2$ or $3 \pi/2$. Since $f(r \cos(\theta), r \sin(\theta))$ does not depend on $r$, this upper bound will not be exceeded either when $r$ escapes to $+\infty$ or when $r$ approaches 0. We have therefore established that
	\[
	\max_{(x, y) \in \big(\mathds{R} \setminus \{0\}\big) \times \mathds{R}} \dfrac{(2 (\alpha + 1) x + y) (x + y)^{2 \alpha} - y y^{2 \alpha}}{(x + y)(x + y)^{2 \alpha} - y y^{2 \alpha}} = 2 \left[ \dfrac{(\alpha + 1) m_\alpha + 1}{m_\alpha + 1} \right].
	\]
	This maximum provides the desired optimal $N$ in the inequality \eqref{eqn:mcpv0}.
\end{proof}

It seems that the Grushin structures behave in such a way that points $q_0$ and $q$ lying on the same horizontal line (with $x_0 \neq 0$) provide the sharpest $N$ where $\beta_t(q_0,q) \geqslant t^N$ holds for all $t \in \interval{0}{1}$. This is also what happens when $\alpha = 1$ (see \cite[Proposition 62.]{subriemint} and \cite[Theorem 8.]{rizzicut} for Grushin half-planes).

We thus expect that the optimal $N$ obtained in \cref{prop:distcoefhor} is sharp. We are able to verify this intuition for singular points, i.e. when $q_0 = (0, y_0)$.

\begin{proposition}
	Let $q_0 = (x_0, y_0) \in \mathds{G}_\alpha$ with $x_0 = 0$ and $q \notin \mathrm{Cut}(q_0)$. The inequality
	\[
	\beta_t(q, q_0) \geqslant t^N
	\]
	holds for all $t \in \interval{0}{1}$ and every $N \geqslant 2 \left[ \dfrac{(\alpha + 1) m_\alpha + 1}{m_\alpha + 1} \right]$.
\end{proposition}

\begin{proof}
	We firstly observe that
	\begin{equation}
		\label{usefulbound}
		N_\alpha := 2 \left[ \dfrac{(\alpha + 1) m_\alpha + 1}{m_\alpha + 1} \right] \geqslant 2 (\alpha + 1),
	\end{equation}
	since $\alpha \geqslant 1 > 0$.
	
	If $x_0 = 0$ and $v_0 = 0$, the formula of the distortion coefficients in \cref{prop:v0} and \cref{usefulbound} yield $\beta_t(q_0,q) = t^{2 (\alpha + 1)} \geqslant t^{N_\alpha}$.
	
	Assume now that $x_0 = 0$ and $v_0 \neq 0$, i.e. $\phi = 0$ or $\phi = \pi_\alpha$, the Jacobian determinant \eqref{jacobianexp} is given by
	\begin{equation}
	\label{detjacx0}
	[x, y](t) = t \dfrac{A^{2 \alpha}}{\alpha \omega} \Bigl[\sin_\alpha(\omega t) - \omega t \cos_\alpha(\omega t) \Bigr].
	\end{equation}
	It follows from \eqref{detjacx0} that
	\begin{equation}
	\label{betax0}
	\beta_t(q, q_0) = t \dfrac{g(\omega t)}{g(\omega)}
	\end{equation}
	where we have set $g(z) := \sin_\alpha(z) - z \cos_\alpha(z)$.
	We first note that $g(0) = 0$. Then, we compute
	\begin{equation*}
	g'(z) = \alpha z\sin_\alpha^{2 (\alpha - 1)}(z) \sin_\alpha(z)
	\end{equation*}
	to find that $g'(z) > 0$ for every $ z \in \ointerval{0}{\pi_\alpha}$ and $\alpha \geqslant 1$. Therefore, the functions $g$ is strictly increasing and positive.
	
	We want to prove that \eqref{betax0} is greater than $t^{N_\alpha}$. In the same way as in the proof of \cref{prop:distcoefhor}, we know that the desired inequality holds if and only if we have
	\[
	G(z) := (N_\alpha - 1) g(z) - z g'(z) \geqslant 0 \text{ for all } z \in \interval{0}{\pi_\alpha}. 
	\]
	We can see that $G(0) = 0$, and
	\begin{equation*}
	G'(z) = \alpha z \sin_\alpha^{2 (\alpha - 1)}(z) \left[ (N - 3) \sin_\alpha(z) - (2 \alpha - 1) z \cos_\alpha(z) \right].
	\end{equation*}
	From \cref{usefulbound}, we deduce that
	\begin{equation*}
		G'(z) \geqslant \alpha (2 \alpha - 1) z \sin_\alpha^{2 (\alpha - 1)}(z) \left[ \sin_\alpha(z) - z \cos_\alpha(z) \right].
	\end{equation*}
	Therefore, $G'(z)$ is non-negative and so is $G(z)$, for all $\in \interval{0}{\pi_\alpha}$.
\end{proof}

By analysing in more detail and looking at the graph of the distortion coefficients \eqref{distAphiom}, it would indeed seem to us that the relevant condition is also satisfied when $x_0 \neq 0$. We therefore propose \cref{conjMCP}. A proof of this could require further work, potentially involving a more comprehensive study of the $(2, 2\alpha)$-trigonometric functions.					

\nocite{*}

\printbibliography


\end{document}